\documentclass[11pt]{article}
\usepackage{authblk}
\usepackage[utf8]{inputenc}
\usepackage[margin=1in]{geometry}
\usepackage{amssymb, amsmath, amsthm}
\usepackage{mathtools}
\usepackage{graphicx}
\usepackage{enumitem}
\usepackage{tikz}
\usepackage{mathrsfs}
\usepackage{comment}
\usepackage[T1]{fontenc}
\usepackage{hyperref}

\theoremstyle{plain}
\newtheorem{lemma}{Lemma}[section]
\newtheorem{claim}{Claim}[section]
\newtheorem{trick}{Fact}[section]
\newtheorem{theorem}[lemma]{Theorem}
\newtheorem{proposition}[lemma]{Proposition}
\newtheorem{corollary}[lemma]{Corollary}
\newtheorem{definition}[lemma]{Definition}
\newtheorem{theorem-definition}[lemma]{Theorem-Definition}

\theoremstyle{remark}
\newtheorem{remark}{Remark}
\newtheorem{example}{Example}

\newcommand{\Z}{\mathbb Z}

\newcommand{\R}{\mathbb R}
\newcommand{\C}{\mathbb C}
\newcommand{\E}{\mathbb E}
\newcommand{\m}{\mathbf m}
\renewcommand{\P}{\mathcal P}
\newcommand{\NC}{\mathcal{NC}}
\newcommand{\Luk}{\mathbf L}
\newcommand{\even}{\textrm{even}}
\newcommand{\odd}{\textrm{odd}}
\renewcommand{\a}{\mathbf a}
\renewcommand{\b}{\mathbf b}
\newcommand{\cc}{\mathbf c}
\newcommand{\x}{\mathbf x}

\newcommand{\q}{\mathbf q}
\renewcommand{\bar}{\,|\,}

\title{Cumulants in rectangular finite free probability and beta-deformed singular values}

\date{\today}

\author{Cesar Cuenca\thanks{cesar.a.cuenk@gmail.com}}
\affil{Department of Mathematics, The Ohio State University, Columbus, OH, USA.}

\newtheoremstyle{named}{}{}{\itshape}{}{\bfseries}{.}{.5em}{\thmnote{#3}#1}
\theoremstyle{named}
\newtheorem*{namedtheorem}{}

\begin{document}

\maketitle

\abstract{
Motivated by the $(q,\gamma)$-cumulants, introduced by Xu~\cite{X_2023} to study $\beta$-deformed singular values of random matrices, we define the $(n,d)$-rectangular cumulants for polynomials of degree~$d$ and prove several moment-cumulant formulas by elementary algebraic manipulations; the proof naturally leads to quantum analogues of the formulas.
We further show that the $(n,d)$-rectangular cumulants linearize the $(n,d)$-rectangular convolution from Finite Free Probability and that they converge to the $q$-rectangular free cumulants from Free Probability in the regime where $d\to\infty$, $1+n/d\to q\in[1,\infty)$.
As an application, we employ our formulas to study limits of symmetric empirical root distributions of sequences of polynomials with nonnegative roots.
One of our results is akin to a theorem of Kabluchko~\cite{K_2022} and shows that applying the operator $\exp(-\frac{s^2}{n}x^{-n}D_xx^{n+1}D_x)$, where $s>0$, asymptotically amounts to taking the rectangular free convolution with the rectangular Gaussian distribution of variance $qs^2/(q-1)$.
}


\section{Introduction}

\subsection{Preface}

The \emph{asymptotic freeness}, due by Voiculescu~\cite{V_1991}, was the first bridge discovered between Free Probability and Random Matrix Theory. This principle roughly states that certain large independent Hermitian random matrices behave like free random variables.
In particular, the empirical spectral distribution of sums of independent Hermitian random matrices tends to the free convolution of measures, as the sizes of the matrices grow to infinity.

As discovered in~\cite{MSS_2022}, asymptotic freeness can already be perceived when considering random matrices of a fixed size $d\times d$, and this has given rise to the growing field of \emph{Finite Free Probability}, see e.g.~\cite{AP_2018,M_2021,G_2022,G_2024,GM_2022,AGVP_2023}, etc., for different aspects of the theory.
The main idea in the subject is to take the expected characteristic polynomials; then the sum of random matrices leads to a binary operation on polynomials of degree~$d$, called the \emph{symmetric additive convolution}, which in some sense tends to the free convolution of measures, as $d\to\infty$.
Thus, Finite Free Probability can be understood as a refinement of Free Probability, which includes the additional positive integer parameter $d$, namely the degree of the polynomials.
Of special interest for us will be the combinatorial viewpoint on Finite Free Probability, suggested by~\cite{AP_2018}, where \emph{finite free cumulants} of polynomials were introduced and shown to linearize the symmetric additive convolution; such cumulant-based viewpoint was coined for Free Probability in~\cite{S_1994}.

From an entirely different point of view, a novel one-parameter deformation of free convolution\footnote{This new convolution takes as input two probability measures and outputs a signed measure in general; however, it is conjectured that the output is always a probability measure.} and associated cumulants were discovered in the setting of high temperature $\beta$-ensembles, by considering random $\beta$-sums of $d\times d$ matrices~\cite{BGCG_2022}.
The high temperature regime is when $d\to\infty$ and $\beta d/2\to\gamma$ simultaneously; the parameter $\gamma\in (0,\infty)$ is the one that remains present in the new version of free convolution.
It turns out, unexpectedly, that the convolution and cumulants in this $\gamma$-setting agree with the symmetric additive convolution and finite free cumulants under the formal identification of parameters $\gamma\leftrightarrow -d$.
The present paper is an effort to find applications of this remarkable coincidence.
Specifically, motivated by the theory of high temperature random \emph{$\beta$-sums of rectangular matrices}~\cite{X_2023}, we define and study cumulants for the operation of \emph{rectangular convolution of polynomials}~\cite{GM_2022} from Finite Free Probability.
We hope that the link considered sheds new light on both subjects.

Finally, it is worth mentioning that there are other recent works~\cite{CDM_2023, CD_2025a, CD_2025b} that study the quantized setting of $\beta$-tensor products of symmetric group representations and which have led to yet different and new one-parameter deformations of free convolution and cumulants.
Whether those objects are also linked to Finite Free Probability is an open problem.

\subsection{General setting}

Let $t,u$ be generic real numbers.
Let $\a=(a_2,a_4,\dots)$ be a real sequence, indexed by the even positive integers.
Define its associated sequence of \emph{cumulants} $\kappa^\a=(\kappa^\a_2,\kappa^\a_4,\dots)$ and \emph{moments} $m^\a=(m^\a_2,m^\a_4,\dots)$ by means of the following generating series identities:
\begin{equation}\label{def_cm}
\begin{aligned}
\exp\left( \sum_{\ell=1}^\infty{ \frac{\kappa^\a_{2\ell}}{\ell} z^{2\ell} } \right)
&= 1 + \sum_{n=1}^\infty{ \frac{a_{2n}}{(t)_n(u)_n}z^{2n} },\\
\exp\left( t\sum_{k=1}^\infty{ \frac{m^\a_{2k}}{k} z^{2k} } \right) &= 1 + \sum_{n=1}^\infty{ a_{2n}z^{2n} }.
\end{aligned}
\end{equation}
We employed above the notation~$(v)_n := \prod_{i=0}^{n-1}{(v+i)}$.
Our first main result compiles combinatorial formulas expressing any $m_{2k}^\a$ in terms of $\kappa^\a_2,\kappa^\a_4,\cdots$.

\begin{namedtheorem}[Main Result~I]\emph{(See Thms.~\ref{thm:m_k}, \ref{thm:m_k_2} and \ref{thm:m_k_3} in the text.)}
For any $k\in\Z_{\ge 1}$, the moment $m^\a_{2k}$ is the constant term of the formal power series
\begin{equation}\label{eqn_xu_0}
\left( \partial_{t,u} + *_g \right)^{2k-1}g(z),
\end{equation}
where $g(z) := \sum_{\ell=1}^\infty{\kappa^\a_{2\ell} z^{2\ell-1}}$; the operator $*_g$ is multiplication by $g(z)$; and $\partial_{t,u}$ is the operator defined on monomials by $\partial_{t,u} 1:= 0$, and $\partial_{t,u} z^{2m-1} := (u+m-1)z^{2m-2}$, $\partial_{t,u} z^{2m} := (t+m)z^{2m-1}$, for all $m\in\Z_{\ge 1}$.

Alternatively, $m_{2k}^\a$ can be expressed as a sum of $\kappa^\a$-weighted \textbf{odd \L{}ukasiewicz paths} of length $2k$, and as a sum of $\kappa^\a$-weighted \textbf{even set partitions} of $\{1,2,\dots,2k\}$.
The definitions of these objects are in Sec.~\ref{sec:prelims} and the precise theorem statements are in Sec.~\ref{sec:cums_moms}.
\end{namedtheorem}

We point out that the formula in Eqn.~\eqref{eqn_xu_0} can be deduced from the results of~\cite{X_2023}.
However, the proof there involves sophisticated tools, such as the BC-type Dunkl operators and transform.
It was an open problem to find a direct combinatorial proof and this is our first contribution!
Indeed, our proofs in Sec.~\ref{sec:cums_moms} are elementary algebraic manipulations starting from~\eqref{def_cm}.
In Sec.~\ref{sec:q_analogues}, we prove a $\q$-analogue of Main Result~I that involves the nontrivial notion of $\q$-composition~\cite{G_1982,J_1996}; see Thm.~\ref{thm:q_analogue}.
It is worth noting that this $\q$-version was obtained by following our combinatorial proof, but it seems difficult to deduce it from the connection to BC-type Dunkl theory.
But now that this $\q$-analogue is proved, a natural question arises: is it related to some version of BC-type Dunkl theory powered by Koornwinder symmetric polynomials~\cite{R_2005}, or to some $\q$-deformation of free probability~\cite{KS_2009,LT_2016}?

The motivation for the definition~\eqref{def_cm} of cumulants and moments is only apparent when one considers the operation of convolution of sequences.
Let $\a=(a_2,a_4,\dots)$, $\b=(b_2,b_4,\dots)$ and $\cc=(c_2,c_4,\dots)$ be real sequences satisfying
\begin{equation}\label{eqn:convolution_intro}
\frac{c_{2n}}{(t)_n(u)_n} = \sum_{k=0}^n\frac{a_{2k}}{(t)_k(u)_k}\cdot\frac{b_{2n-2k}}{(t)_{n-k}(u)_{n-k}},
\,\text{ for all }n\in\Z_{\ge 1},
\end{equation}
where $a_0:=1$, $b_0:=1$.
The formula~\eqref{eqn:convolution_intro} is the convolution of sequences $(\a,\b)\mapsto\cc$.
Upon certain specialization of $t,u$, it turns out that this map is exactly the \emph{rectangular convolution} of polynomials from finite free probability~\cite{GM_2022}, as shown in Sec.~\ref{sec:rectangular_conv}.
We use~\eqref{def_cm} as motivation to define the \emph{rectangular (finite free) cumulants} that linearize the rectangular convolution.
This contribution, explained in Sec.~\ref{subsec:example_1}, should be compared to the paper of Arizmendi-Perales~\cite{AP_2018} that works out similar ideas for the \emph{symmetric additive convolution}.
On the other hand, the map $(\a,\b)\mapsto\cc$ defined by~\eqref{eqn:convolution_intro} is equivalent to the corresponding map of moments $(m^\a,m^\b)\mapsto m^\cc$, which under different specializations of $t,u$, describes precisely the rectangular free convolution of compactly supported symmetric probability measures (in terms of moment sequences) discovered by Benaych-Georges~\cite{BG_2009}, as well as the $(q, \gamma)$-convolution of Xu~\cite{X_2023}; this connection is explained in Appendix~\ref{subsec:appendix_1}.
Our first main result provides moment-cumulant formulas applicable to both subjects.

\subsection{Special setting: Rectangular convolution}\label{subsec:example_1}

For any $n\in\Z_{\ge 0}$ and monic polynomials $p(x) = x^{d} + \sum_{i=1}^d{a_{2i}\,x^{d-i}}$, $r(x) = x^{d} + \sum_{i=1}^d{b_{2i}\,x^{d-i}}$ of degree~$d$, the \emph{$(n,d)$-rectangular convolution} of $p(x)$ and $r(x)$ is defined by\footnote{The map $\boxplus^n_d$ from~\cite{GM_2022} is bilinear, as are all other related maps from Finite Free Probability~\cite{MSS_2022}. W.l.o.g.~in this paper, we only consider monic polynomials as inputs for these operations.}
\begin{equation*}
(p\boxplus_d^n r)(x) := x^{d} + \sum_{k=1}^d{ x^{d-k}\sum_{i+j=k}{ \frac{(d-i)!(d-j)!}{d!(d-k)!} \frac{(n+d-i)!(n+d-j)!}{(n+d)!(n+d-k)!}\, a_{2i} b_{2j} } }.
\end{equation*}
Alternatively, the $(n,d)$-rectangular convolution can be defined as follows.
If $D:=\frac{d}{dx}$, there are unique polynomials $P(x)$, $R(x)$ of degree $d$ such that $p(x)=P(x^{-n}Dx^{n+1}D)x^d$ and $r(x)=R(x^{-n}Dx^{n+1}D)x^d$; then $(p\boxplus^n_dr)(x) = P(x^{-n}Dx^{n+1}D)R(x^{-n}Dx^{n+1}D)x^d$. The case $n=0$ is an observation of Mirabelli~\cite{Mir_2021}, also stated in~\cite{MSS_2022}.
This crucial fact will be used in the proof of our Main Result~III on the asymptotic theory of nonnegative real-rooted polynomials and, in fact, it was essential in the conception of its statement.

This convolution is motivated by random matrix theory.
Specifically, if we consider independent $d\times (n+d)$ random matrices $A, B$ with deterministic singular values and $p(x)$, $r(x)$ are the characteristic polynomials of $AA^*, BB^*$, then $(p\boxplus_d^nr)(x)$ is the expected characteristic polynomial of $(A+UBV)(A+UBV)^*$, where $U\in U(d)$ and $V\in U(n+d)$ are Haar-distributed; see~\cite{GM_2022,MSS_2022} for details.

A last related result that we will mention is that if all roots of $p(x), r(x)$ are real and nonnegative, then the same is true of $(p\boxplus_d^n r)(x)$; see~\cite[Thm.~2.3]{GM_2022}.

If $p(x)$ is a polynomial of degree~$d$ with nonnegative real roots of the form $\alpha_1^2,\dots,\alpha_d^2$, for some $\alpha_1,\dots, \alpha_d\ge 0$, then its \emph{symmetric empirical root distribution} is defined by
\begin{equation}\label{eq:symmetric_distribution}
\widetilde\mu[p] := \frac{1}{2d}\sum_{i=1}^d{(\delta_{\alpha_i} + \delta_{-\alpha_i})}.
\end{equation}
Equivalently, this can be thought of as the empirical root distribution of the polynomial $\widetilde{p}(x) := p(x^2)$ with roots $\pm\alpha_1,\dots,\pm\alpha_d$.

For any polynomial $p(x) = x^{d} + \sum_{i=1}^d{a_{2i}\,x^{d-i}}$ of degree~$d$, define its \emph{$(n,d)$-rectangular (finite free) cumulants} $K^{n,d}_2[p],K^{n,d}_4[p],\dots,K^{n,d}_{2d}[p]$ by
\begin{equation}\label{eq:def_rectangular_cumulants}
K^{n,d}_{2\ell}[p] := \ell\cdot[z^{2\ell}]\,\ln\left( 1 + \sum_{i=1}^d{\frac{a_{2i}}{(-d)_i(-d-n)_i}z^{2i}} \right),
\quad\text{for all $\ell=1,\dots,d$},
\end{equation}
where $[z^{2\ell}]f(z)$ is the coefficient of $z^{2\ell}$ in the formal power series $f(z)\in\R[[z]]$.
This definition is motivated by Eqn.~\eqref{def_cm}, upon the specialization of parameters
\begin{equation*}
t=-d,\qquad u=-d-n.
\end{equation*}

\begin{namedtheorem}[Main Result~II]
The $(n,d)$-rectangular cumulants ``linearize'' the $(n,d)$-rectangular convolution and ``finitize'' the $q$-rectangular free cumulants, in the following sense:

(a) For any two monic polynomials $p(x)$ and $r(x)$ of degree~$d$, we have
\begin{equation*}
K^{n,d}_{2\ell}\big[ p\boxplus^n_d r \big] = K^{n,d}_{2\ell}[p] + K^{n,d}_{2\ell}[r],\quad\text{for all }\ell=1,\dots,d.
\end{equation*}

(b) Let $q\in [1,\infty)$, let $\{p_d(x)\}_{d\ge 1}$ be such that each $p_d(x)$ is a monic polynomial of degree~$d$ with nonnegative real roots, and let $\mu$ be a symmetric\footnote{A measure $\mu$ on $\R$ is said to be \emph{symmetric} if $\mu(A)=\mu(-A)$, for all Lebesgue sets $A\subseteq\R$, where $-A:=\{-a\mid a\in A\}$.} probability measure on the real line with finite moments of all orders.
Then $\widetilde\mu[p_d]\to\mu$, as $d\to\infty$, in the sense of moments, if and only if
\begin{equation*}
\lim_{d\to\infty,\, 1+\frac{n}{d}\to q}{ (-d)^{2\ell-1}K^{n,d}_{2\ell}[p_d] } = q^{-\ell}\kappa^q_{2\ell}[\mu],
\quad\textrm{for all $\ell\in\Z_{\ge 1}$}.
\end{equation*}

(c) Let $q\in[1,\infty)$ and let $\mu,\nu$ be compactly supported symmetric probability measures on the real line. Let $\{p_d(x)\}_{d\ge 1}$, $\{r_d(x)\}_{d\ge 1}$ be such that $p_d(x)$, $r_d(x)$ are monic polynomials of degree~$d$ with all their roots being real and nonnegative.
Assume that $\widetilde\mu[p_d]\to\mu$, $\widetilde\mu[r_d]\to\nu$, as $d\to\infty$, in the sense of moments.
Then $\widetilde\mu\big[ p_d\boxplus^n_d r_d \big]\to\mu\boxplus_q\nu$, in the sense of moments, in the regime where $n,d\to\infty$ and $1+\frac{n}{d}\to q$.
\end{namedtheorem}

In part~(b), $\kappa^q_{2\ell}[\mu]$ denote the $q$-rectangular free cumulants of $\mu$, whereas in part~(c), $\mu\boxplus_q\nu$ denotes the $q$-rectangular free convolution of $\mu$ and $\nu$; both of these objects were defined and studied in~\cite{BG_2009}, but see also Sec.~\ref{sec:bg_section} for our notations and a summary of the facts needed.

Parts (a), (b) and (c) are Thm.~\ref{thm:linearity}, Thm.~\ref{thm:limit_moments_cumulants} and Cor.~\ref{cor:d_infinity} in the text, respectively.
This second main result should be compared to the analogous statements~\cite[Prop.~3.6 \& Cor.~5.5]{AP_2018} for \emph{finite free cumulants}.
That paper was motivated by the \emph{$d$-finite $R$-transform}~\cite{M_2021}.
We note that~\cite[Def.~3.7]{G_2024} is the definition of a \emph{rectangular finite $R$-transform}, denoted in that paper by $\mathcal{R}_{\mathbb{S}p}^{d,\lambda}(s)$. Our rectangular cumulants~\eqref{eq:def_rectangular_cumulants} are the coefficients of this transform, up to a constant prefactor and identification of parameters $\lambda=d/(n+d)$.
Then part (a) also follows from~\cite[Thm.~3.9]{G_2024}, which shows that the rectangular finite $R$-transform linearizes the rectangular convolution of polynomials.
Likewise, part (b) follows from the coefficient-wise limit~\cite[Thm.~6.1]{G_2024}.
Our proofs are independent and employ only our combinatorial formulas from the first main result.

\subsection{Application to the asymptotic theory of polynomials}

Let $\{p_d(x)\}_{d\ge 1}$ be a sequence of monic polynomials, such that each $p_d(x)$ is of degree~$d$ and its roots are of the form $(\alpha_1^{(d)})^2,\dots, (\alpha_d^{(d)})^2$, for some nonnegative real $\alpha_1^{(d)},\dots,\alpha_d^{(d)}\ge 0$. Let
\begin{equation*}
\widetilde\mu[p_d] := \frac{1}{d}\sum_{i=1}^d{ \left( \delta_{\alpha_i^{(d)}} + \delta_{-\alpha_i^{(d)}} \right) },\qquad d\in\Z_{\ge 1},
\end{equation*}
be the corresponding symmetric empirical root distributions.
The following result should be compared to \cite[Thm.~2.11]{K_2022}; that theorem relates to symmetric additive convolution in the same way that our following next result relates to rectangular convolution.

\begin{namedtheorem}[Main Result~III]
Assume the following convergence in the sense of moments
\begin{equation*}
\widetilde\mu[p_d] \to \mu,\ \text{ as }d\to\infty,
\end{equation*}
where $\mu$ is a compactly supported symmetric probability measure on the real line.
Consider arbitrary $s\in(0,\infty)$, $q\in(1,\infty)$, $n\in\Z_{\ge 0}$, and denote the derivative by $D:=\frac{d}{dx}$. Then
\begin{equation}\label{eq:conclusion_3}
\widetilde\mu\left[ \exp\left(-\frac{s^2}{n} x^{-n}Dx^{n+1}D\right)p_d(x) \right] \longrightarrow \mu\boxplus_q\lambda^{(q)}_{qs^2/(q-1)},\quad\text{as }\,n,d\to\infty,\ 1+\frac{n}{d}\to q,
\end{equation}
in the sense of moments. Above, $\lambda^{(q)}_{qs^2/(q-1)}$ denotes the probability measure on $\R$ with density
\begin{equation*}
\frac{d\lambda^{(q)}_{qs^2/(q-1)}(x)}{dx} = \frac{\sqrt{4qs^4 - \big((q-1)x^2-(q+1)s^2\big)^2}}{2\pi s^2|x|}\cdot\mathbf{1}_{A(q,s)}(x),
\end{equation*}
where $\mathbf{1}_{A(q,s)}(x)$ denotes the indicator function of the set
\begin{equation*}
A(q,s) := \left[ -\frac{\sqrt{q-1}}{\sqrt{q}-1}s,\, -\frac{\sqrt{q-1}}{\sqrt{q}+1}s \right]
\cup\left[ \frac{\sqrt{q-1}}{\sqrt{q}+1}s,\, \frac{\sqrt{q-1}}{\sqrt{q}-1}s \right]\subseteq\R.
\end{equation*}
\end{namedtheorem}

The measure $\lambda^{(q)}_{qs^2/(q-1)}$ is the $q$-rectangular analogue of the centered Gaussian distribution of variance $qs^2/(q-1)$, discussed in~\cite{BG_2007,BG_2009}.
The proof is in Sec.~\ref{sec:asymptotic_real_roots}.

\subsection{Organization of the paper}

Besides this introduction, the present paper has five other sections plus an appendix. In Sec.~\ref{sec:prelims}, we introduce the combinatorial background of set partitions, \L{}ukasiewicz paths and the $q$-rectangular free convolution.
In Sec.~\ref{sec:cums_moms} we find moment-cumulant formulas for sequences related to each other by means of Eqn.~\eqref{def_cm}; in particular, Main Result~I is proved in this section.
In Sec.~\ref{sec:cumulants_finite_free}, we begin with the applications to finite free probability, namely we define the rectangular cumulants and prove part (a) from Main Result~II.
Then in Sec.~\ref{sec:applications}, we study limits of empirical root distributions of polynomials with nonnegative real roots and prove Main Result~III, as well as parts (b)-(c) from Main Result~II.
In Sec.~\ref{sec:q_analogues}, we consider a $\q$-generalization to Main Result~I.
Finally, in Appendix~\ref{appendix} we briefly explain the versions of cumulants related to high temperature $\beta$-deformed singular values and eigenvalues.

\subsection*{A remark on the notation}

The reader might be curious as to why most sequences in this paper are labeled by the positive \emph{even} integers (notably $(m_2,m_4,\dots)$ and $(\kappa_2,\kappa_4,\dots)$), as opposed to all positive integers.
This boils down to the fact that~\cite{BG_2009} and~\cite{X_2023} chose to study the symmetric empirical root distributions~\eqref{eq:symmetric_distribution} of polynomials with nonnegative real roots, as opposed to the empirical root distributions.
Such choice implies that all measures that appear in their studies have vanishing odd moments $m_1=m_3=\cdots=0$ and vanishing odd cumulants $\kappa_1=\kappa_3=\cdots=0$.
We decided to stick to the notations used by these authors, thus leading to the current presentation.
If one considered instead the empirical root distributions, we could rewrite part of the theory for probability measures supported on $[0,\infty)$ instead; for them, generally, odd-indexed $m_{2n-1}$ and $\kappa_{2n-1}$ will be nonzero.

\subsection*{Acknowledgments}

The author would like to thank Jiaming Xu, Hoi Nguyen, Volodymyr Chub and Stanislav Surmylo for helpful conversations, as well as two anonymous referees for helpful comments.
This article is a follow-up to the final report by Chub--Surmylo of a project from the research program Yulia's Dream; the author is grateful to the organizers for allowing him to play the role of research mentor.
This work was partially supported by the NSF grant DMS-2348139.

\section{Preliminaries}\label{sec:prelims}

\subsection{Set partitions}\label{sec:set_partitions}

Let $n\in\Z_{\ge 1}$. A \emph{set partition} of $[n]:=\{1,2,\cdots,n\}$ is an unordered collection of pairwise disjoint nonempty subsets $B_1,\cdots, B_k$ of $[n]$ such that $[n] = B_1\cup\cdots\cup B_k$; the corresponding set partition is denoted $\pi=\{B_1, \cdots, B_k\}$. The subsets $B_1,\cdots,B_k$ are called the \emph{blocks} of $\pi$. The cardinalities of the blocks are denoted $|B_1|, \cdots, |B_k|$. Also, denote the number of blocks of $\pi$ by $\#(\pi)$; in our running example, $\#(\pi)=k$.
Finally, the set of all set partitions of $[n]$ will be denoted by $\P(n)$.

For any $\pi,\sigma\in\P(n)$, we write that $\pi\ge\sigma$, or $\sigma\le\pi$, if $\sigma$ is a \emph{refinement} of $\pi$, i.e.~if any block of $\sigma$ is contained in some block of $\pi$; this is called the \emph{reverse refinement order} of set partitions.

We say that $\pi=\{B_1,\dots,B_k\}\in\P(n)$ is a \emph{noncrossing set partition} if whenever $1\le a<b<c<d\le n$, $\,a, c\in B_i$, and $b,d\in B_j$, then $i=j$. The set of all noncrossing set partitions of $[n]$ will be denoted by $\NC(n)$.

When $n=2n'$ is even, we distinguish yet another class of set partitions.
We say that $\pi\in\P(2n')$ is \emph{even} if all its blocks have even cardinalities.
Denote the set of all even set partitions of $[2n']$ by $\P^\even(2n')\subseteq\P(2n')$.
Also, denote the the set of all even noncrossing set partitions of $[2n']$ by $\NC^\even(2n') := \P^\even(2n')\cap\NC(2n')$.

\subsection{Generating series and set partitions}

For any sequence of variables $x_1,x_2,\dots$, and any $\pi\in\P(n)$, we denote
\begin{equation*}
x_\pi := \prod_{B\in\pi}{x_{|B|}}.
\end{equation*}
We will make use of the following classical result, see e.g.~\cite[Ch.~1]{S_2012} or \cite[Sec.~1.2]{MS_2017}.

\begin{proposition}\label{xs_ys_prop}
Let $x_1,x_2,\dots$ and $y_1,y_2,\dots$ be two sequences related by
\begin{equation}\label{xs_ys}
1 + \sum_{n=1}^\infty{ x_n \frac{z^n}{n!} } = \exp\left( \sum_{n=1}^\infty{y_n \frac{z^n}{n!}} \right).
\end{equation}
Then for any $n\in\Z_{\ge 1}$:
\begin{align*}
x_n &= \sum_{\pi\in\P(n)}{y_\pi},\\
y_n &= \sum_{\pi\in\P(n)}{ (-1)^{\#(\pi)-1}(\#(\pi) - 1)!\, x_\pi }.
\end{align*}
\end{proposition}

\subsection{\L{}ukasiewicz paths}\label{sec:luk_paths}

By definition, a \emph{\L{}ukasiewicz path of length $n$} is a lattice path in $\Z^2$ from $(0,0)$ to $(n,0)$ that never goes below the $x$-axis and has $n$ steps, each being of the form $(1,j)$, for some $j\in\Z$, $j\ge -1$.
A step $(1,j)$, for some $j\ge 1$ (resp.~$j=0$ or $j=-1$) is called an up step (resp.~horizontal or down step); also, if $(i,m)$ is the vertex of some \L{}ukasiewicz path, we say that this vertex is at height~$m$.
We denote by $\Luk(n)$ the set of all \L{}ukasiewicz paths of length $n$.

\begin{proposition}[Prop.~9.8 from~\cite{NS_2006}]\label{prop:bijection}
Given $\pi=\{B_1,\dots,B_k\}\in\NC(n)$, let $a_i\in[n]$ be the smallest element in $B_i$, for $1\le i\le k$.
Next, construct a lattice path $\Lambda(\pi)$ that begins at $(0,0)$ and is followed by $n$ steps: its $m$-th step is $(1,|B_i|-1)$, if $m=a_i$ for some $i$, while it is $(1,-1)$ otherwise.
Then the constructed lattice path $\Lambda(\pi)$ is a \L{}ukasiewicz path of length $n$ and the map $\pi\mapsto\Lambda(\pi)$ is a bijection between $\NC(n)$ and $\Luk(n)$.
\end{proposition}

We say that a \L{}ukasiewicz path is \emph{odd} if all its steps are of the form $(1,2k-1)$, for some $k\in\Z_{\ge 0}$. Note that any odd \L{}ukasiewicz path must have even length. If $n=2n'$, then denote by $\Luk^\odd(2n')\subseteq\Luk(2n')$ the set of odd \L{}ukasiewicz paths of length $2n'$.
For example, some of the odd \L{}ukasiewicz paths belonging to $\Luk^\odd(6)$ are shown in Fig.~\ref{fig:luk} (ignore the step-labels in the figure, for now).

A quick analysis of the bijective map $\NC(n)\to\Luk(n)$ described in Prop.~\ref{prop:bijection} in the case when $n=2n'$ shows that the image of $\NC^\even(2n')$ is exactly $\Luk^\odd(2n')$.
Let us record this observation, as well as others that will be useful, in the next corollary.

\begin{corollary}\label{cor:bijection_2}
The map $\pi\mapsto\Lambda(\pi)$ described in Prop.~\ref{prop:bijection} furnishes a bijection between $\NC^\even(2n')$ and $\Luk^\odd(2n')$. Moreover, we have the following correspondences:

\smallskip

$\bullet$ If the $m$-th step of $\Lambda(\pi)\in\Luk^\odd(2n')$ is an up step, then $m$ is the smallest element of some block of $\pi$.

$\bullet$ If the $m$-th step of $\Lambda(\pi)\in\Luk^\odd(2n')$ is a down step from some odd height $2s+1$ to the even height $2s$, then $m$ is even and is not the smallest element of any block of $\pi$.

$\bullet$ If the $m$-th step of $\Lambda(\pi)\in\Luk^\odd(2n')$ is a down step from some even height $2s$ to the odd height $2s-1$, then $m$ is odd and is not the smallest element of any block of $\pi$.
\end{corollary}

\subsection{q-rectangular free convolution and q-rectangular cumulants}\label{sec:bg_section}

We will present here the cumulant-based combinatorial side of the theory of $q$-rectangular free probability; the material here and further details, including proofs to the statements made, are in~\cite{BG_2009}.\footnote{The parameter $\lambda\in[0,1]$ that is used in~\cite{BG_2009} is related to our parameter $q\in[1,\infty)$ by the relation $\lambda=q^{-1}$.}
Let $q\in [1,\infty)$ be a real number and let $\mu$ be a symmetric probability measure on $\R$ with finite moments of all orders.
Since $\mu$ is symmetric, it has vanishing odd moments; its even moments, on the other hand, will be denoted by
\begin{equation*}
m_{2k}[\mu] := \int_\R{ x^{2k}\mu(dx) },\quad k\in\Z_{\ge 1}.
\end{equation*}
Further, define the \emph{$q$-rectangular free cumulants} of $\mu$, to be denoted by $\kappa^q_{2\ell}[\mu]$, $\ell\in\Z_{\ge 1}$, recursively by the equations
\begin{equation}\label{eq:q_rectangular_cumulants}
m_{2k}[\mu] = \sum_{\pi\in\NC^\even(2k)}{ q^{-\even(\pi)}\prod_{B\in\pi}{\kappa^q_{|B|}[\mu]} },\quad k\in\Z_{\ge 1},
\end{equation}
where $\even(\pi)$ denotes the number of blocks of $\pi$ whose smallest element is an even number in $[2k]$.
For example, Eqns.~\eqref{eq:q_rectangular_cumulants} for $k=1,2,3$ lead to the first three $q$-rectangular free cumulants:
\begin{gather*}
\kappa^q_2[\mu] = m_2[\mu],\qquad
\kappa^q_4[\mu] = m_4[\mu] - (1+q^{-1})\cdot m_2[\mu]^2,\\
\kappa^q_6[\mu] = m_6[\mu] - 3(1+q^{-1})\cdot m_4[\mu]m_2[\mu] + (2+3q^{-1}+2q^{-2})\cdot m_2[\mu]^3.
\end{gather*}

\begin{theorem-definition}\label{thm_def}
Let $\mu,\nu$ be two compactly supported symmetric probability measures on $\R$, with $q$-rectangular free cumulants denoted by $\kappa_{2\ell}^q[\mu],\kappa_{2\ell}^q[\nu]$, for all $\ell\in\Z_{\ge 1}$.
Then there exists a unique probability measure, to be denoted $\mu\boxplus_q\nu$, with finite moments of all orders and corresponding $q$-rectangular free cumulants being
\begin{equation*}
\kappa_{2\ell}^q \big[ \mu\boxplus_q\nu \big] = \kappa_{2\ell}^q[\mu] + \kappa_{2\ell}^q[\mu],\quad\textrm{for all }\ell\in\Z_{\ge 1}.
\end{equation*}
Moreover, $\mu\boxplus_q\nu$ is compactly supported and symmetric.
The probability measure $\mu\boxplus_q\nu$ is called the \textbf{$q$-rectangular free convolution} of $\mu$ and $\nu$.
\end{theorem-definition}

There exists a symmetric probability measure, uniquely determined by its moments and with only its second $q$-rectangular free cumulant being nonzero.
This is the $q$-rectangular analogue of the centered Gaussian distribution for the theory of $q$-rectangular free probability.
It will be denoted by $\lambda^{(q)}_{\sigma^2}$, if it has variance $\sigma^2>0$.
The density of $\lambda^{(q)}_{\sigma^2}$, stated next, was obtained in~\cite[Thm.~4.3]{BG_2007} for $\sigma^2=1$; see also~\cite[Sec.~3.10.2]{BG_2009} for general $\sigma^2>0$.\footnote{There is a small typo in the formula from~\cite[Sec.~3.10.2]{BG_2009}.}

\begin{lemma}\label{lem:helpful_measure}
Let $\sigma\in (0,\infty)$, $q\in[1,\infty)$ be arbitrary. The measure $\lambda^{(q)}_{\sigma^2}$ on $\R$ with density
\begin{equation*}
\frac{d\lambda^{(q)}_{\sigma^2}(x)}{dx} = \frac{\sqrt{4q\sigma^4 - (qx^2 - (q+1)\sigma^2)^2}}{2\pi\sigma^2|x|}\cdot
\mathbf{1}_{B(q,\sigma)}(x),
\end{equation*}
where
\begin{equation*}
B(q,\sigma) := \left[ -\big(1+q^{-\frac{1}{2}}\big)\sigma,\, -\big(1-q^{-\frac{1}{2}}\big)\sigma \right]\cup\left[ \big(1-q^{-\frac{1}{2}}\big)\sigma,\, \big(1+q^{-\frac{1}{2}}\big)\sigma \right]\subseteq\R,
\end{equation*}
is a compactly supported symmetric probability measure, uniquely determined by its $q$-rectangular free cumulants $\displaystyle\kappa^q_{2\ell}\Big[ \lambda^{(q)}_{\sigma^2}\Big] = \delta_{\ell, 1}\cdot\sigma^2$.
\end{lemma}

\section{Moment-cumulant formulas}\label{sec:cums_moms}

Hereinafter we will use the Pochhammer symbol, defined by $(v)_n:=\prod_{i=0}^{n-1}(v+i)$, for all $n\in\Z_{\ge 0}$.

In this section, we consider $\a=(a_2,a_4,\dots)$, $\pmb{\kappa}=(\kappa_2,\kappa_4,\dots)$ and $\m=(m_2,m_4,\dots)$, indexed by the positive even integers, which we call the sequences of coefficients, cumulants and moments, respectively, and which are related to each other by:
\begin{align}
\exp\left( \sum_{\ell=1}^\infty{ \frac{\kappa_{2\ell}}{\ell} z^{2\ell} } \right)
&= 1 + \sum_{n=1}^\infty{ \frac{a_{2n}}{(t)_n(u)_n}z^{2n} },\label{eqn:a_k}\\
\exp\left( t\sum_{k=1}^\infty{ \frac{m_{2k}}{k} z^{2k} } \right) &= 1 + \sum_{n=1}^\infty{ a_{2n}z^{2n} }.\label{eqn:a_m}
\end{align}
In the equations above, $t, u$ can be treated as formal parameters or as generic real numbers.
The interest on these equations lies on their significance in the high temperature limits of $\beta$-deformed singular values, as explained in Appendix~\ref{appendix}.
They will also be used in the upcoming two sections, in connection to finite free probability.

Below, we shall prove several transition formulas between $\pmb{\kappa}$ and $\m$.
There are three such results, namely Theorems~\ref{thm:m_k}, \ref{thm:m_k_2} and \ref{thm:m_k_3}.
It turns out that the moments $m_{2k}$ can be expressed as polynomials on the cumulants $\kappa_2,\kappa_4,\dots,\kappa_{2k}$, with coefficients being polynomials in $\Z_{\ge 0}[t, u]$. This is why we can prove more combinatorial formulas for the transition $\pmb{\kappa}\mapsto\mathbf{m}$, while only the less explicit Thm.~\ref{thm:m_k} gives a formula for the transition $\mathbf{m}\mapsto\pmb{\kappa}$.

\begin{example}
From~\eqref{eqn:a_k}--\eqref{eqn:a_m}, one finds that the first few moments in terms of cumulants are
\begin{align*}
m_2 =&\ u\kappa_2,\\
m_4 =&\ u(t+1)(u+1)\kappa_4 + u(t+u+1)\kappa_2^2,\\
m_6 =&\ u(t+1)(u+1)(t+2)(u+2)\kappa_6 + 3u(t+1)(u+1)(t+u+2)\kappa_2\kappa_4\\
&+ u(t^2+u^2+3tu+3t+3u+2)\kappa_2^3,
\end{align*}
while the first few cumulants in terms of moments are
\begin{align*}
\kappa_2 =&\ \frac{m_2}{u},\\
\kappa_4 =&\ \frac{m_4}{u(t+1)(u+1)} - \frac{t+u+1}{u^2(t+1)(u+1)}\,m_2^2,\\
\kappa_6 =&\ \frac{m_6}{u(t+1)(u+1)(t+2)(u+2)} - \frac{3(t+u+2)}{u^2(t+1)(u+1)(t+2)(u+2)}\,m_2m_4\\
&+ \frac{2t^2+2u^2+3tu+6t+6u+4}{u^3(t+1)(u+1)(t+2)(u+2)}\,m_2^3.
\end{align*}
\end{example}

\subsection{In terms of even set partitions}

Eqn.~\eqref{eqn:a_k} defines each $a_{2n}$ as a function of $\kappa_2,\dots,\kappa_{2n}$ and each $\kappa_{2n}$ as a function of $a_2,\dots,a_{2n}$. More precisely:

\begin{lemma}\label{lem:a_k}
For any $n\in\Z_{\ge 1}$, we have
\begin{align*}
a_{2n} &= \frac{(t)_n(u)_n}{(2n)!} \sum_{\pi\in\P^\even(2n)} \left( 2^{\#(\pi)} \prod_{B\in\pi}{(|B|-1)!} \right)\cdot\kappa_\pi,\\
\kappa_{2n} &= \frac{1}{2(2n-1)!} \sum_{\pi\in\P^\even(2n)} \left( (-1)^{\#(\pi)-1} (\#(\pi)-1)! 
\prod_{B\in\pi} \frac{|B|!}{(t)_{\frac{|B|}{2}}(u)_{\frac{|B|}{2}}} \right)\cdot a_\pi.
\end{align*}
\end{lemma}

\begin{proof}
Our Eqn.~\eqref{eqn:a_k} is identical to Eqn.~\eqref{xs_ys} under the following variable identification:
\begin{gather*}
y_{2n-1}=0, \qquad y_{2n}=2(2n-1)!\,\kappa_{2n}, \quad\text{for all } n\in\Z_{\ge 1},\\
x_{2n-1}=0, \qquad x_{2n}=\frac{(2n)!}{(t)_{n}(u)_n}\, a_{2n}, \quad\text{for all } n\in\Z_{\ge 1}.
\end{gather*}
As a result, Prop.~\ref{xs_ys_prop} gives the following formulas:
\begin{align*}
a_{2n} = \frac{(t)_n(u)_n}{(2n)!}\, x_{2n}
&= \frac{(t)_n(u)_n}{(2n)!} \sum_{\pi\in\P(2n)}\prod_{B\in\pi} y_{|B|} \\
&= \frac{(t)_n(u)_n}{(2n)!} \sum_{\pi\in\P^\even(2n)} \prod_{B\in\pi} y_{|B|} \\
&= \frac{(t)_n(u)_n}{(2n)!} \sum_{\pi\in\P^\even(2n)} 2^{\#(\pi)} \prod_{B\in\pi}{(|B|-1)!\cdot\kappa_\pi},
\end{align*}
and this proves the first equality from the statement of the proposition.
Observe that the equality between lines 1 and 2 above follows from the fact that $y_{2k-1}=0$, for all $k\in\Z_{\ge 1}$, so the set partitions with some block of odd cardinality gives a zero contribution and therefore the sum over all set partitions $\P(2n)$ can in fact be restricted to the set $\P^\even(2n)$ of even set partitions.
The second equality is proved analogously, by using the second identity from Prop.~\ref{xs_ys_prop}.
\end{proof}

Likewise, Eqn.~\eqref{eqn:a_m} defines each $a_{2n}$ as a function of $m_2,\dots,m_{2n}$ and each $m_{2n}$ as a function of $a_2,\dots,a_{2n}$.
Indeed, the following is proved similarly to the previous proposition; it is also proved in~\cite[Thm.~4.1.1]{CS_2024}.

\begin{lemma}\label{lem:m_a}
For any $n\in\Z_{\ge 1}$, we have
\begin{align*}
a_{2n} &= \frac{1}{(2n)!}\sum_{\pi\in\P^\even(2n)}{ \left( (2t)^{\#(\pi)}\prod_{B\in\pi}{(|B|-1)!} \right)\cdot m_\pi},\\
m_{2n} &= \frac{1}{2(2n-1)!\cdot t}\,\sum_{\pi\in\P^\even(2n)}{ \left( (-1)^{\#(\pi)-1}(\#(\pi)-1)! \prod_{B\in\pi}{|B|!} \right)\cdot a_\pi}.
\end{align*}
\end{lemma}

As a consequence of Lemmas~\ref{lem:a_k} and~\ref{lem:m_a}, we obtain the following.

\begin{theorem}\label{thm:m_k}
For any $n\in\Z_{\ge 1}$, we have
\begin{multline*}
m_{2n} = \frac{1}{2(2n-1)!\cdot t}\ 
\sum_{\sigma\in\P^\even(2n)} \bigg( 2^{\#(\sigma)}\prod_{V\in\sigma}{(|V|-1)!}\\
\cdot\sum_{\pi\colon\pi\ge\sigma}{ (-1)^{\#(\pi)-1}(\#(\pi)-1)! \prod_{B\in\pi}{(t)_{\frac{|B|}{2}}(u)_{\frac{|B|}{2}}} }
\bigg)\cdot\kappa_\sigma,
\end{multline*}
and
\begin{equation*}
\kappa_{2n} = \frac{1}{2(2n-1)!}\sum_{\sigma\in\P^\even(2n)}
\bigg( (2t)^{\#(\sigma)}\prod_{V\in\sigma}(|V|-1)!
\sum_{\pi\colon\pi\ge\sigma}{\frac{(-1)^{\#(\pi)-1}(\#(\pi)-1)!}{\prod_{B\in\pi}(t)_{\frac{|B|}{2}}(u)_{\frac{|B|}{2}}}} \bigg) \cdot m_\sigma.
\end{equation*}
\end{theorem}

This theorem is analogous to~\cite[Thm.~4.2]{AP_2018} (see also~\cite[Thm.~4.2.1]{CS_2024}) and the proof is very similar, so we will omit it.
Let us only remark that the inner summations in both formulas are over set partitions $\pi$ such that $\sigma$ is a refinement of $\pi$, and since $\sigma\in\P^\even(2n)$, then automatically $\pi\in\P^\even(2n)$, as well; in particular, if $B\in\pi$, then $\frac{|B|}{2}\in\Z_{\ge 1}$.

\subsection{In terms of operators applied to formal power series}

\begin{definition}\label{def:partial_tu}
Let $\partial_{t,u}$ be the linear operator on the space $\R[[z]]$ of formal power series on $z$, defined on monomials by
\begin{align*}
\partial_{t,u} z^{2m+1} &:= (u+m)z^{2m},\ \text{ for all $m\in\Z_{\ge 0}$},\\
\partial_{t,u} z^{2m} &:= (t+m)z^{2m-1},\ \text{ for all $m\in\Z_{\ge 1}$},\\
\partial_{t,u} 1 &:= 0.
\end{align*}
\end{definition}

\begin{theorem}\label{thm:m_k_2}
For any $k\in\Z_{\ge 1}$, we have
\begin{equation}\label{eqn_xu}
m_{2k} = [z^0] \left( \partial_{t,u} + *_g \right)^{2k-1}g(z),
\end{equation}
where $\partial_{t,u}$ is the operator from Definition~\ref{def:partial_tu},
\begin{equation*}
g(z) := \sum_{\ell=1}^\infty{\kappa_{2\ell} z^{2\ell-1}},
\end{equation*}
the operator $*_g$ is multiplication by $g(z)$, and $[z^0]\colon\R[[z]]\to\R$ picks up the constant term of a power series, i.e.~$[z^0]f(z) = f(0)$.
\end{theorem}

This result follows from \cite[Thms.~4.8 and~5.8]{X_2023}.
In fact, that paper defines for certain symmetric probability measures with even moments $m_2,m_4,\dots$ the \emph{$(q,\gamma)$-cumulants} by means of Eqn.~\eqref{eqn_xu}, with $t=\gamma$ and $u=q\gamma$, for some $\gamma\in(0,\infty)$, $q\in[1,\infty)$.
Then it is proven that the moments and $(q,\gamma)$-cumulants are related to each other by~\eqref{eqn:a_k}--\eqref{eqn:a_m}, and alternatively by~\eqref{eqn_xu}.
As the class of probability measures for which this equivalence holds is large enough, the desired equivalence of identities must hold for any sequences $\mathbf{m}$ and $\pmb{\kappa}$, regardless of whether they come from probability measures.
The disadvantage of the approach in \cite{X_2023} is that it relies heavily on the Dunkl transform and operators associated to the root system of type BC, so it is a conceptually advanced and technical proof.
Here, we offer a simple proof that involves only the manipulation of power series.
The same technique also furnishes a proof of the equivalence between Thms.~3.10 and~3.11 in \cite{BGCG_2022} that was proved there by the use of the Dunkl transform and Dunkl operators of type A.

\begin{proof}[Proof of Theorem~\ref{thm:m_k_2}]
\textbf{Step 1 (Formula for $a_{2n}$ in terms of $\kappa_{2\ell}$'s and the operator $\partial_{t,u}$).}
Let
\begin{equation}\label{G_def}
G(z):=\sum_{\ell\ge 1}{ \frac{\kappa_{2\ell}}{\ell} z^{2\ell} },
\end{equation}
so that $\frac{G'(z)}{2}=g(z)$.
By taking derivatives to both sides of~\eqref{eqn:a_k}, we deduce
\begin{equation}\label{eqn:a_k_derivative}
g(z)e^{G(z)} = \sum_{n=1}^\infty{ \frac{na_{2n}}{(t)_n(u)_n}z^{2n-1} }.
\end{equation}
The key observation here are the equations
\begin{align*}
\frac{\partial_{t,u}z^{2n-1}}{(t)_n(u)_n} &= \frac{z^{2n-2}}{(t)_n(u)_{n-1}},\text{ for all }n\ge 1,\\
\frac{\partial_{t,u}z^{2n-2}}{(t)_n(u)_{n-1}} &= \frac{z^{2n-3}}{(t)_{n-1}(u)_{n-1}},\text{ for all }n\ge 2,
\end{align*}
which are derived from the definition of the operator $\partial_{t,u}$.
Thus, by applying $\partial_{t,u}^{2n-1}$ to \eqref{eqn:a_k_derivative} and then taking the constant term on both sides, we obtain
\begin{equation}\label{eqn:step1}
a_{2n} = \frac{t}{n}\cdot [z^0]\,\partial_{t,u}^{2n-1}(g(z)\cdot e^{G(z)}), \text{ for all }n\ge 1.
\end{equation}

\textbf{Step 2 (Recursive relation between $a_{2n}$'s and $m_{2k}$'s).}
By taking derivatives with respect to $z$ from both sides in~\eqref{eqn:a_m}, we have
\begin{equation*}
t\cdot\sum_{k=1}^\infty{ m_{2k} z^{2k-1} }\cdot\exp\left( t\sum_{k=1}^\infty{ \frac{m_{2k}}{k} z^{2k} } \right)
= \sum_{n=1}^\infty{ na_{2n}z^{2n-1} }.
\end{equation*}
Then by replacing the exponential in this equation with the right hand side of \eqref{eqn:a_m}, we have
\begin{equation*}
t\cdot\sum_{k=1}^\infty{ m_{2k} z^{2k-1} }\cdot\left( 1+\sum_{n=1}^\infty{ a_{2n}z^{2n} } \right)
= \sum_{n=1}^\infty{ na_{2n}z^{2n-1} }.
\end{equation*}
By comparing the coefficients of $z^{2n-1}$ on both sides, we deduce
\begin{equation}\label{eqn:step2}
\frac{na_{2n}}{t} = m_{2n} + \sum_{k=1}^{n-1}{ m_{2k}a_{2n-2k} }.
\end{equation}
From this equation, note that $m_{2n}$ is a function of $m_{2n-2}, \dots, m_2,a_{2n},\dots,a_2$.
We can then prove the desired
\begin{equation}\label{eqn:step2_2}
m_{2k} = [z^0] \left( \partial_{t,u} + *_g \right)^{2k-1}(g(z)), \text{ for all }k\ge 1,
\end{equation}
by induction on $k$. The base case $k=1$ is trivial, so only the inductive step remains.
For the latter, observe that it suffices to verify that the $a_{2n}$'s defined by~\eqref{eqn:step1} and the $m_{2k}$'s defined by~\eqref{eqn:step2_2} satisfy the recurrence relation \eqref{eqn:step2}. This is achieved in Step 4, after some preparations.

\smallskip
\textbf{Step 3 (The auxiliary operator $d_{t,u}$).} We can write
\begin{equation}\label{eqn:d_tu}
\partial_{t,u} = d_{t,u} + \frac{1}{2}\frac{d}{dz},
\end{equation}
where $d_{t,u}$ is the linear operator on $\R[[z]]$ defined by
\begin{align*}
d_{t,u} z^{2m+1} &:= \left(u-\frac{1}{2}\right)\cdot z^{2m},\quad\text{for all $m\in\Z_{\ge 0}$},\\
d_{t,u} z^{2m} &:= t\cdot z^{2m-1},\quad\text{for all $m\in\Z_{\ge 1}$},\\
d_{t,u} 1 &:= 0.
\end{align*}

\begin{claim}\label{claim:step_3}
For $G(z)$ as in Eqn.~\eqref{G_def} and any formal power series $h(z)\in\R[[z]]$, we have
\begin{equation}\label{eqn:step3}
\partial_{t,u}\big( h(z)e^{G(z)} \big) = (\partial_{t,u} + *_g)(h(z))\cdot e^{G(z)} + [z^0]h(z)\cdot d_{t,u}\big( e^{G(z)} \big).
\end{equation}
\end{claim}
\begin{proof}[Proof of Claim~\ref{claim:step_3}]
Observe that $d_{t,u}$ acts as follows on certain products of powers of $z$:
\begin{align*}
d_{t,u}(z^{2a+1}z^{2b}) &= d_{t,u}(z^{2a+2b+1}) = \left(u-\frac{1}{2}\right)\cdot z^{2a+2b} = (d_{t,u}z^{2a+1})\cdot z^{2b},\,\text{ for all $a,b\in\Z_{\ge 0}$},\\
d_{t,u}(z^{2a}z^{2b}) &= d_{t,u}(z^{2a+2b}) = t\cdot z^{2a+2b-1} = (d_{t,u}z^{2a})\cdot z^{2b},\,\text{ for all $a\in\Z_{\ge 1},b\in\Z_{\ge 0}$},\\
d_{t,u}(1\cdot z^{2b}) &= \mathbf{1}_{\{b\ne 0\}}\cdot tz^{2b-1} = 1\cdot d_{t,u}z^{2b},\,\text{ for all $b\in\Z_{\ge 0}$}.
\end{align*}
By linearity of the operator $d_{t,u}$, the previous equations imply
\begin{equation}\label{eqn:step3_1}
d_{t,u}(h(z)u(z^2)) = (d_{t,u}h(z))\cdot u(z^2) + [z^0]h(z)\cdot d_{t,u}(u(z^2)),
\end{equation}
where $u(z^2)$ represents a power series with only even powers of $z$.
Since $G(z)=\sum_{\ell\ge 1}{ \frac{\kappa_{2\ell}}{\ell} z^{2\ell} }$ only has even powers of $z$, then so does $e^{G(z)}$, and we can replace $u(z^2)$ in~\eqref{eqn:step3_1} by $e^{G(z)}$ to obtain
\begin{equation}\label{eqn:step3_2}
d_{t,u}\big( h(z)e^{G(z)} \big) = (d_{t,u}h(z))\cdot e^{G(z)} + [z^0]h(z)\cdot d_{t,u}\big( e^{G(z)} \big).
\end{equation}

On the other hand, by the Leibniz rule and the fact that $\frac{G'(z)}{2}=g(z)$, we have
\begin{equation}\label{eqn:step3_3}
\frac{1}{2}\cdot \frac{d}{dz}\left( h(z)e^{G(z)} \right) = \left( \frac{1}{2}\cdot \frac{d}{dz}(h(z)) + g(z)\cdot h(z) \right)\cdot e^{G(z)}.
\end{equation}
Adding \eqref{eqn:step3_2} and \eqref{eqn:step3_3} proves the desired \eqref{eqn:step3} and ends the proof of the claim.
\end{proof}

\textbf{Step 4 (Conclusion of the argument).}
As mentioned in the last sentence of Step 2, to finish the proof of the theorem, it will suffice to show that the $a_{2n}$'s defined by~\eqref{eqn:step1} and $m_{2k}$'s defined by~\eqref{eqn:step2_2} satisfy the recursive relations~\eqref{eqn:step2}.
Note that $g(z)=\sum_{\ell\ge 1}{\kappa_{2\ell}z^{2\ell-1}}$ only has odd powers of $z$, so $[z^0](g(z))=0$.
Then it follows from Claim~\ref{claim:step_3} that
\begin{equation}\label{eqn:step4_1}
\partial_{t,u}\big( g(z)e^{G(z)} \big) = (\partial_{t,u} + *_g)(g(z))\cdot e^{G(z)}.
\end{equation}
We want to take the constant term to both sides of this equation.
For the right hand side, we use that $[z^0](u(z)v(z))=[z^0](u(z))\cdot [z^0](v(z))$, for any $u,v\in\R[[z]]$, and also that $[z^0]e^{G(z)}=1$, because $G(z)$ does not have constant terms.
Then~\eqref{eqn:step4_1} yields
\begin{equation*}
[z^0]\,\partial_{t,u}\big( g(z)e^{G(z)} \big) = [z^0](\partial_{t,u} + *_g)(g(z)).
\end{equation*}
From \eqref{eqn:step1} and \eqref{eqn:step2_2}, this gives $\frac{a_2}{t} = m_2$, which is exactly the desired Eqn.~\eqref{eqn:step2} for $n=1$.

Next, apply $\partial_{t,u}$ to~\eqref{eqn:step4_1}. From Claim~\ref{claim:step_3}, this gives
\begin{equation}\label{eqn:step4_2}
\partial_{t,u}^2\big( g(z)e^{G(z)} \big) = (\partial_{t,u} + *_g)^2(g(z))\cdot e^{G(z)}
+ [z^0](\partial_{t,u} + *_g)(g(z))\cdot d_{t,u}\big( e^{G(z)} \big)
\end{equation}
Now note that $\partial_{t,u}$ decreases the degree of a power of $z$ by $1$. Also, because $g(z)$ has only odd powers of $z$, it follows that $(\partial_{t,u}+*_g)(g(z))$ only has even powers of $z$, and $(\partial_{t,u}+*_g)^2(g(z))$ only has odd powers of $z$; in particular, $[z^0]\,(\partial_{t,u}+*_g)^2(g(z))=0$.
Now we can apply $\partial_{t,u}$ to~\eqref{eqn:step4_2} and use Claim~\ref{claim:step_3} again to get
\begin{equation}\label{eqn:step4_3}
\partial_{t,u}^3\big( g(z)e^{G(z)} \big) = (\partial_{t,u} + *_g)^3(g(z))\cdot e^{G(z)}
+ [z^0](\partial_{t,u} + *_g)(g(z))\cdot\partial_{t,u}d_{t,u}\big( e^{G(z)} \big).
\end{equation}
Proceeding inductively in the same fashion, we obtain
\begin{multline}\label{eqn:step4_general}
\partial_{t,u}^{2n-1}\big( g(z)e^{G(z)} \big) = (\partial_{t,u} + *_g)^{2n-1}(g(z))\cdot e^{G(z)}\\
+ \sum_{m=1}^{n-1}{ [z^0](\partial_{t,u} + *_g)^{2m-1}(g(z))\cdot\partial_{t,u}^{2n-2m-1}d_{t,u}\big( e^{G(z)} \big) },
\end{multline}
for all $n\in\Z_{\ge 1}$. Taking the constant terms from both sides yields
\begin{multline}\label{eqn:step4_general_2}
[z^0]\,\partial_{t,u}^{2n-1}(g(z)e^{G(z)}) = [z^0]\,(\partial_{t,u} + *_g)^{2n-1}(g(z))\\
+ \sum_{m=1}^{n-1}{ [z^0](\partial_{t,u} + *_g)^{2m-1}(g(z))\cdot [z^0]\,\partial_{t,u}^{2n-2m-1}d_{t,u}(e^{G(z)}) }.
\end{multline}
To simplify $[z^0]\,\partial_{t,u}^{2n-2m-1}d_{t,u}(e^{G(z)})$, we use the following obvious fact.

\begin{trick}\label{trick_step_3}
Let $k\in\Z_{\ge 1}$ and let $\partial_1,\dots,\partial_k$ be operators on $\R[[z]]$, each of which is a multiple of $\frac{d}{dz}$.
For any $f(z)=\sum_{n=0}^\infty{a_nz^n}\in\R[[z]]$, if we denote $[z^k]f(z):=a_k$, then we have
\begin{equation*}
[z^0]\,\partial_1\cdots\partial_k(f(z)) = \partial_1\cdots\partial_k\Big( [z^k]f(z) \cdot z^k \Big).
\end{equation*}
\end{trick}
Since $\partial_{t,u}$ and $d_{t,u}$ are multiples of $\frac{d}{dz}$, we can then deduce
\begin{align}
[z^0]\,\partial_{t,u}^{2n-2m-1}d_{t,u}(e^{G(z)}) &= \partial_{t,u}^{2n-2m-1}d_{t,u}\Big( [z^{2n-2m}]e^{G(z)} \cdot z^{2n-2m} \Big)\nonumber\\
&= \partial_{t,u}^{2n-2m-1}\left( \frac{t}{2n-2m}\cdot\frac{d}{dz} \right)\Big( [z^{2n-2m}]e^{G(z)} \cdot z^{2n-2m} \Big)\nonumber\\
&= \frac{t}{2n-2m}\cdot [z^0]\,\partial_{t,u}^{2n-2m-1}\frac{d}{dz}\big( e^{G(z)} \big)\nonumber\\
&= \frac{t}{n-m}\cdot [z^0]\,\partial_{t,u}^{2n-2m-1}\big( g(z)e^{G(z)} \big),\label{eqn:after_trick}
\end{align}
where the first and third equalities follow from Fact~\ref{trick_step_3}, the second equality holds because $d_{t,u}z^{2\ell}=\frac{t}{2\ell}\cdot\frac{d}{dz}(z^{2\ell})$, by the definition of the action of $d_{t,u}$ and the last equality is because $\frac{G'(z)}{2}=g(z)$.
Plugging~\eqref{eqn:after_trick} into~\eqref{eqn:step4_general_2} gives:
\begin{multline}\label{eqn:step4_general_3}
[z^0]\,\partial_{t,u}^{2n-1}(g(z)e^{G(z)}) = [z^0]\,(\partial_{t,u} + *_g)^{2n-1}(g(z))\\
+ \sum_{m=1}^{n-1}{ [z^0](\partial_{t,u} + *_g)^{2m-1}(g(z))\cdot
\left( \frac{t}{n-m}\cdot [z^0]\,\partial_{t,u}^{2n-2m-1}(g(z)e^{G(z)}) \right) },
\end{multline}
Finally notice that, by using \eqref{eqn:step1} and \eqref{eqn:step2_2}, this relation precisely matches \eqref{eqn:step2} and, hence, finishes the proof.
\end{proof}

\subsection{In terms of \L{}ukasiewicz paths}

\begin{theorem}\label{thm:m_k_3}
For any $k\in\Z_{\ge 1}$, we have
\begin{multline}\label{eqn:m_k_3}
m_{2k} = \sum_{P\in\Luk^\odd(2k)} \prod_{s\ge 0}{(u+s)^{\#\textrm{down steps of $P$ from height }(2s+1)}} \\
\cdot\prod_{s\ge 1}{(t+s)^{\#\textrm{down steps of $P$ from height }(2s)}}
\prod_{s\ge 1}{\kappa_{2s}^{\#\textrm{up steps }(1,2s-1)\textrm{ of $P$}}}.
\end{multline}
In other words, for any odd \L{}ukasiewicz path $P$ of length $2k$, associate to each of its steps a weight according to the following rules. To each up step $(1,2s-1)$, associate the weight $\kappa_{2s}$ and to each down step $(1,-1)$ from height $(j+1)$ to height $j$ the weight $(u+s)$, if $j=2s$ is even, and the weight $(t+s)$, if $j=2s-1$ is odd.\footnote{By definition, there are no horizontal steps in an odd \L{}ukasiewicz path.}
Finally, associate to $P$ the product of all weights of its steps and denote it by $w_{t, u}(P)$.
Then $m_{2k}$ is the sum of weights $w_{t,u}(P)$, as $P$ ranges over all odd \L{}ukasiewicz paths of length $2k$.
\end{theorem}

This result is equivalent to \cite[Thm.~5.5]{X_2023}. Indeed, that result gives an equivalent formula in terms of noncrossing even set partitions, which relate to odd \L{}ukasiewicz paths as explained by Cor.~\ref{cor:bijection_2}.
The advantage of our theorem is that the formula in terms of \L{}ukasiewicz paths is more natural and also bears resemblance to universal formulas in the related discrete context of high temperature beta-partitions~\cite{CDM_2023,CD_2025a,CD_2025b}.

\begin{proof}[Proof of Theorem~\ref{thm:m_k_3}]
Begin with Thm.~\ref{thm:m_k_2}:
\begin{equation*}
m_{2k} = [z^0] \left( \partial_{t,u} + *_g \right)^{2k-1}g(z).
\end{equation*}
We can write
\begin{equation*}
\left( \partial_{t,u} + *_g \right)^{2k-1}g(z) = \left( \partial_{t,u} + *_g \right)^{2k-1}*_g(1)
= \sum_{a_1,\dots,a_{2k}}{a_{2k}\cdots a_2a_1(1)},
\end{equation*}
where the sum is over ``words'' $a_{2k}\cdots a_2a_1$ of length $2k$ and each $a_i$ is either the operator $\partial_{t,u}$ or $*_g$, but the last ``letter'' is $a_1=*_g$.
Since $g(z)=\sum_{\ell\ge 1}{\kappa_{2\ell}z^{2\ell-1}}$, this sum can be further refined and written as:
\begin{equation}\label{sum_expression}
\left( \partial_{t,u} + *_g \right)^{2k-1}g(z) = \sum_{b_1,\dots,b_{2k}}{b_{2k}\cdots b_2b_1(1)},
\end{equation}
where each $b_i$ is either $\partial_{t,u}$ or $*_{\kappa_{2\ell}z^{2\ell-1}}$ (operator of multiplication by $\kappa_{2\ell}z^{2\ell-1}$), for some $\ell\in\Z_{\ge 1}$, and $b_1$ is of the latter type.

Let us define the degree of $b_i$ to be $\deg(b_i):=-1$, if $b_i=\partial_{t,u}$, and $\deg(b_i):=2\ell-1$, if $b_i=*_{\kappa_{2\ell}z^{2\ell-1}}$.
For any specific term $b_{2k}\cdots b_2b_1(1)$, if some ``suffix'' $b_j\cdots b_1$ has $\deg(b_j)+\cdots+\deg(b_1)<0$, there is a smallest such suffix (with smallest $j$), then necessarily $b_j=\partial_{t,u}$ and $b_{j-1}\cdots b_1(1)$ is a constant, therefore $b_j\cdots b_1(1)=0$ and also $b_{2k}\cdots b_2b_1(1)=0$.
As a result, the sum~\eqref{sum_expression} can be restricted to words $b_{2k}\cdots b_2b_1$ with the additional constraint that all suffixes $b_j\cdots b_1$ have $\deg(b_j)+\cdots+\deg(b_1)\ge 0$.
Note that if $\deg(b_{2k})+\cdots+\deg(b_1)=n\ge 0$, then $b_{2k}\cdots b_1(1)$ is a constant multiple of $z^n$.
Thus, by taking the constant term of \eqref{sum_expression}, we obtain
\begin{equation}\label{sum_expression_2}
m_{2k} = [z^0]\left( \partial_{t,u} + *_g \right)^{2k-1}g(z) = \sum_{c_1,\dots,c_{2k}}{c_{2k}\cdots c_2c_1(1)},
\end{equation}
where the sum is now over words $c_{2k}\cdots c_1$ of length $2k$, satisfying:

\smallskip

$\bullet$ each $c_i$ is either $\partial_{t,u}$ or $*_{\kappa_{2\ell}z^{2\ell-1}}$, for some $\ell\in\Z_{\ge 1}$;

$\bullet$ $c_1=*_{\kappa_{2\ell}z^{2\ell-1}}$, for some $\ell\in\Z_{\ge 1}$;

$\bullet$ all suffixes $c_j\cdots c_1$ have $\deg(c_j)+\cdots+\deg(c_1)\ge 0$;

$\bullet$ $\deg(c_{2k})+\cdots+\deg(c_1)=0$.

\smallskip

\noindent Associate to each term $c_{2k}\cdots c_2c_1(1)$ in the sum~\eqref{sum_expression_2} an odd \L{}ukasiewicz path of length $2k$,
\begin{equation*}
P = \left( w_0=(0,0)\to w_1\to w_2\to\cdots\to w_{2k}=(2k,0) \right),
\end{equation*}
with certain edge-weights according to the following procedure. If $c_1=*_{\kappa_{2\ell}z^{2\ell-1}}$, then set $w_1 := (1,2\ell-1)$ and put an edge-weight of $\kappa_{2\ell}$ on $(w_0\to w_1)$. In general, if $w_1,\dots, w_{j-1}$ have already been chosen and $c_j=*_{\kappa_{2\ell}z^{2\ell-1}}$, then set $w_j:=w_{j-1}+(1,2\ell-1)$ and assign an edge-weight of $\kappa_{2\ell}$ on $(w_{j-1}\to w_j)$.
On the other hand, if $c_j=\partial_{t,u}$, then set $w_j:=w_{j-1}+(1,-1)$. The weight of the edge $(w_{j-1}\to w_j)$ in this case will depend on the parity of $\deg(c_{j-1})+\cdots+\deg(c_1)$: if this is even and equal to $2s$, assign the weight $(t+s)$ to the edge, whereas if it is odd and equal to $2s+1$, assign the weight $(u+s)$ to the edge.
By induction, one shows that for $j=1,2,\dots,2k$, the suffix $c_j\cdots c_1(1)$ equals $r_jz^{s_j}$, where $s_j=\deg(c_1)+\cdots+\deg(c_j)$ and $r_j$ is the multiplication of edge-weights of $(w_0\to w_1),\dots,(w_{j-1}\to w_j)$.
In particular, the term $c_{2k}\cdots c_1(1)$ in~\eqref{sum_expression_2} equals the weight of the whole constructed \L{}ukasiewicz path; this ends the proof.
\end{proof}

\begin{figure}
\begin{center}
\includegraphics[width=0.6\textwidth]{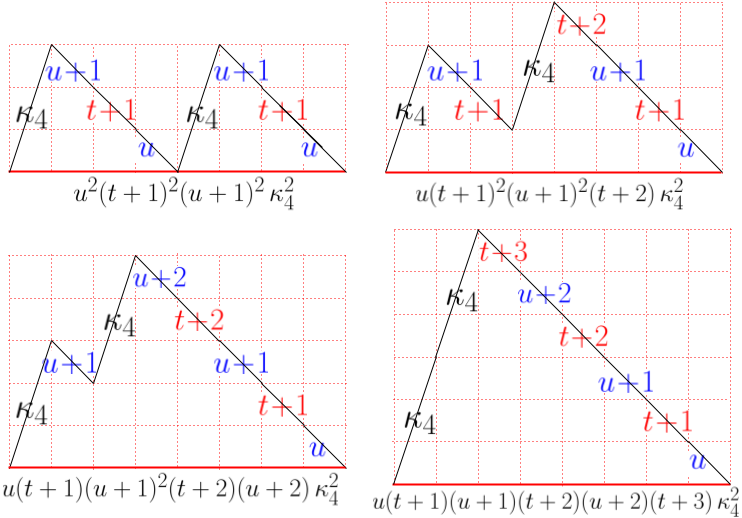}
\caption{The four odd \L{}ukasiewicz paths of length $6$ with two up steps $(1,3)$ and six down steps $(1,-1)$. All edge-weights and weights of the \L{}ukasiewicz paths are also shown.}\label{fig:luk}
\end{center}
\end{figure}

\begin{example}
From the definitions \eqref{eqn:a_k}--\eqref{eqn:a_m}, we find that $m_8$ in terms of cumulants is
\begin{align}
m_8 =&\ u(t+1)(u+1)(t+2)(u+2)(t+3)(u+3)\kappa_8\nonumber\\
&+ 4u(t+1)(u+1)(t+2)(u+2)(t+u+3)\kappa_2\kappa_6\nonumber\\
&+ u(t+1)(u+1)(t+u+3)(2tu+3t+3u+6)\kappa_4^2\label{eqn:coeff}\\
&+ 2u(t+1)(u+1)(3t^2+3u^2+8tu+15t+15u+18)\kappa_2^2\kappa_4\nonumber\\
&+ u(u^3 + t^3 + 6tu^2 + 6t^2u + 6t^2 + 6u^2 + 17tu + 11t + 11u + 6)\kappa_2^4.\nonumber
\end{align}
Let us verify that our Thm.~\ref{thm:m_k_3} correctly predicts the coefficient of $\kappa_4^2$ in $m_8$.
Indeed, a term in the sum~\eqref{eqn:m_k_3} is a multiple of $\kappa_4^2$ if it corresponds to some odd \L{}ukasiewicz path with two up steps $(1,3)$ and six down steps $(1,-1)$. There are exactly four such \L{}ukasiewicz paths; they are exhibited, together with their step-weights, in Fig.~\ref{fig:luk}. The sum of the four weights of the \L{}ukasiewicz paths displayed is indeed $u(t+1)(u+1)(t+u+3)(2tu+3t+3u+6)\kappa_4^2$, matching~\eqref{eqn:coeff}.
\end{example}

\section{Cumulants in finite free probability}\label{sec:cumulants_finite_free}

In the first subsection below, we recall a classical simple result on convolutions of sequences and in the two subsections after, we observe that certain operations from the finite free probability theory can be recast as special cases of convolutions of sequences.

\subsection{Cumulants that linearize convolution of sequences}\label{sec:sequence_cumulants}

For any sequences $\a=(a_1,a_2,\dots)$, $\b=(b_1,b_2,\dots)$, indexed by $\Z_{\ge 1}$, define the \emph{convolution of $\a$ and $\b$} to be the sequence denoted by $\cc=\a*\b$ and defined by
\begin{equation}\label{eqn:sequence_convolution}
\cc=(c_1,c_2,\dots),\qquad c_k := \sum_{i=0}^k{ a_i b_{k-i} },\qquad\text{for all }k\in\Z_{\ge 1},
\end{equation}
where $a_0:=1$, $b_0:=1$.
Further, let $\kappa^\a := (\kappa^\a_1,\kappa^\a_2,\dots)$ be defined from $\a$ by the identity:
\begin{equation}\label{cumulants_coeffs_eq1}
\exp\left( \sum_{\ell=1}^\infty{ \kappa^\a_\ell z^\ell } \right) = 1 + \sum_{n=1}^\infty{a_nz^n}.
\end{equation}
Define $\kappa^\b := (\kappa^\b_1,\kappa^\b_2,\dots)$, $\kappa^\cc := (\kappa^\cc_1,\kappa^\cc_2,\dots)$ similarly.
The following lemma is folklore.

\begin{lemma}\label{thm:basic_convolution}
Let $\a,\b,\cc$ be any sequences and $\kappa^\a,\kappa^\b,\kappa^\cc$ be defined by Eqn.~\eqref{cumulants_coeffs_eq1}. Then
\begin{equation}\label{eqn:linearization}
\cc = \a*\b \Longleftrightarrow \kappa^\cc_\ell = \kappa^\a_\ell + \kappa^\b_\ell, \text{ for all }\ell\in\Z_{\ge 1}.
\end{equation}
\end{lemma}
\begin{proof}
By multiplying the identities
\begin{equation*}
\exp\left( \sum_{\ell=1}^\infty{ \kappa^\a_\ell z^\ell } \right) = 1 + \sum_{n=1}^\infty{a_nz^n},\qquad
\exp\left( \sum_{\ell=1}^\infty{ \kappa^\b_\ell z^\ell } \right) = 1 + \sum_{n=1}^\infty{b_nz^n},
\end{equation*}
we obtain
\begin{equation}\label{conv_eq1}
\exp\left( \sum_{\ell=1}^\infty{ (\kappa^\a_\ell + \kappa^\b_\ell) z^\ell } \right)
= 1 + \sum_{n=1}^\infty{ \left( \sum_{i=0}^n{a_i b_{n-i}} \right) z^n },
\end{equation}
where $a_0,b_0$ are defined to be $a_0:=1$, $b_0:=1$.
The previous identity will be compared to
\begin{equation}\label{conv_eq2}
\exp\left( \sum_{\ell=1}^\infty{ \kappa^\cc_\ell z^\ell } \right) = 1 + \sum_{n=1}^\infty{c_nz^n}.
\end{equation}
Then, by definition, $\cc=\a*\b$ iff the right hand sides of \eqref{conv_eq1} and \eqref{conv_eq2} are equal, iff the left hand sides of \eqref{conv_eq1} and \eqref{conv_eq2} are equal, iff $\kappa^\cc_\ell = \kappa^\a_\ell+\kappa^\b_\ell$, for all $\ell\in\Z_{\ge 1}$.
\end{proof}

\subsection{Symmetric additive convolution}\label{sec:symmetric}

The symmetric additive convolution, defined in~\cite{MSS_2022}, is certain binary operation on real monic polynomials of the same finite degree~$d\in\Z_{\ge 1}$.
Explicitly, if
\begin{equation}\label{polys_p_q_0}
p(x) = x^d+\sum_{i=1}^d{x^{d-i}(-1)^ia_i},\qquad\quad
r(x) = x^d+\sum_{i=1}^d{x^{d-i}(-1)^ib_i},
\end{equation}
then the \emph{symmetric additive convolution} of $p(x)$ and $r(x)$ is defined as
\begin{equation}\label{eq:symmetric_convolution}
p(x) \boxplus_d r(x) := x^d+\sum_{k=1}^d{ x^{d-k}(-1)^k\sum_{i+j=k}{\frac{(d-i)!(d-j)!}{d!(d-k)!}a_ib_j} }.
\end{equation}

On the other hand, to the polynomials $p(x), r(x)$, let us associate the sequences
\begin{equation}\label{sequences_a_b}
\a=\!\left(\frac{a_1}{(-d)_1},\frac{a_2}{(-d)_2},\cdots,\frac{a_d}{(-d)_d},0,0,\cdots\!\right)\!,\quad
\b=\!\left(\frac{b_1}{(-d)_1},\frac{b_2}{(-d)_2},\cdots,\frac{b_d}{(-d)_d},0,0,\cdots\!\right)\!,
\end{equation}
respectively.
Then the convolution of sequences $\cc=\a*\b$, defined according to~\eqref{eqn:sequence_convolution}, is such that all entries of $\cc$ after the $(2d)$ first ones vanish.
Let us write the first $d$ entries of $\cc$ in the form $\frac{c_k}{(-d)_k}$, $1\le k\le d$, for some values of $c_1,\dots,c_d$, i.e:
\begin{equation*}
\cc=\a*\b=\left( \underbrace{\frac{c_1}{(-d)_1},\frac{c_2}{(-d)_2},\dots, \frac{c_d}{(-d)_d}}_{\text{first }d\text{ entries of $a*b$}}, \underbrace{*, *, \dots, *}_{\text{ entries }(d+1),\dots,2d}, \underbrace{0, 0, 0, \dots}_{\text{ entries $>2d$ vanish}} \right).
\end{equation*}
It can be easily verified that the values $c_1,\dots,c_k$ coincide, up to signs, with the coefficients of the symmetric additive convolution $p(x)\boxplus_d r(x)$:
\begin{equation*}
c_k = (-1)^k [x^{d-k}]\left(p(x) \boxplus_d r(x)\right), \text{ for all }k=1,2,\dots, d.
\end{equation*}
Indeed, this follows from plugging the equalities
\begin{equation*}
(-d)_k = (-1)^k\frac{d!}{(d-k)!},\quad (-d)_i = (-1)^i\frac{d!}{(d-i)!},\quad (-d)_j = (-1)^j\frac{d!}{(d-j)!},
\end{equation*}
into \eqref{eqn:sequence_convolution} and comparing with Eqn.~\eqref{eq:symmetric_convolution} that defines symmetric additive convolution.

Motivated by~\eqref{cumulants_coeffs_eq1}, we can define for the monic polynomial $p(x)$ of degree~$d$ in~\eqref{polys_p_q_0} the \emph{finite free cumulants} of $p(x)$ as the first $d$ sequence-cumulants of the sequence $\a$ in \eqref{sequences_a_b}:
\begin{equation}\label{def:finite_cumulants}
K_\ell[p] := \ell\cdot[z^\ell]\,\ln\left( 1 + \sum_{i=1}^d{\frac{a_i}{(-d)_i}z^i} \right),\quad \ell=1,\dots,d,
\end{equation}
where $[z^\ell]f(z)$ denotes the coefficient of $z^\ell$ of any power series $f(z)\in\R[[z]]$.
As a consequence of Lem.~\ref{thm:basic_convolution}, we obtain the following corollary that states that finite free cumulants linearize the symmetric additive convolution.

\begin{corollary}[\cite{AP_2018}]
If $p(x), r(x)$ are monic polynomials of degree~$d$, then
\begin{equation*}
K_\ell\big[ p\boxplus_d r \big] = K_\ell[p] + K_\ell[r],\quad\text{for all $\ell=1,\dots,d$}.
\end{equation*}
\end{corollary}

The quantities $-K_\ell[p]\cdot d^{\ell-1}$, $\ell=1,\dots,d$, coincide with the version of finite free cumulants considered in~\cite{AP_2018}, denoted there by $\kappa^p_{\ell}$.
The previous corollary is therefore equivalent to Prop.~3.6 in that paper, though our proof is different.

\subsection{Rectangular convolution}\label{sec:rectangular_conv}

The rectangular convolution, defined in~\cite{GM_2022}, is certain binary operation on real monic polynomials of the same degree~$d\in\Z_{\ge 1}$. Namely, if
\begin{equation*}
p(x) = x^{d}+\sum_{i=1}^d{x^{d-i}a_{2i}},\qquad\quad
r(x) = x^{d}+\sum_{i=1}^d{x^{d-i}b_{2i}},
\end{equation*}
and $n\in\Z_{\ge 0}$, then the \emph{$(n,d)$-rectangular (additive) convolution} of $p(x)$ and $r(x)$ is defined by
\begin{equation}\label{eq:rectangular_convolution}
p(x)\boxplus_d^n r(x) := x^{d}+\sum_{k=1}^d{ x^{d-k}\sum_{i+j=k}{ \frac{(d-i)!(d-j)!}{d!(d-k)!} \frac{(n+d-i)!(n+d-j)!}{(n+d)!(n+d-k)!}\, a_{2i} b_{2j} } }.
\end{equation}

On the other hand, to the polynomials $p(x)$, $r(x)$, associate the sequences
\begin{align*}
\a &= \left( 0,\frac{a_2}{(-d)_1(-n-d)_1},0,\frac{a_4}{(-d)_2(-n-d)_2},\cdots,\frac{a_{2d}}{(-d)_d(-n-d)_d},0,0,0,\cdots \right),\\
\b &= \left( 0,\frac{b_2}{(-d)_1(-n-d)_1},0,\frac{b_4}{(-d)_2(-n-d)_2},\cdots,\frac{b_{2d}}{(-d)_d(-n-d)_d},0,0,0,\cdots \right).
\end{align*}
Let $\cc=\a*\b$ be the convolution of sequences, defined according to~\eqref{eqn:sequence_convolution}.
Evidently, all odd entries of $\cc$, as well as all entries after the $(4d)$ first ones vanish.
Let us express the first $d$ even entries of $\cc$ as $\frac{c_{2k}}{(-d)_k(-n-d)_k}$, for $k=1,2,\dots,d$, i.e:
\begin{multline*}
\cc = \a*\b = \bigg( \underbrace{0,\frac{c_2}{(-d)_1(-n-d)_1},0,\frac{c_4}{(-d)_2(-n-d)_2},\cdots,0,\frac{c_{2d}}{(-d)_d(-n-d)_d}}_{\text{first $(2d)$ entries of $\a*\b$}},\\
\underbrace{*,*,\dots,*}_{\text{$2d\!+\!1,\dots,4d$}},\ \underbrace{0,0,0,\cdots}_{\text{ entries $>4d$ vanish}} \bigg).
\end{multline*}
Then the values of $c_2,c_4,\dots,c_{2d}$ coincide with the nonzero coefficients of the rectangular convolution $p(x)\boxplus_d^n r(x)$:
\begin{equation*}
c_{2k} = [x^{d-k}]\left(p(x) \boxplus_d^n r(x)\right), \text{ for all }k=1,2,\dots, d.
\end{equation*}
This equality can be verified as in the previous Sec.~\ref{sec:symmetric}.
Motivated again by~\eqref{cumulants_coeffs_eq1}, we introduce the notion of cumulants in this setting.

\begin{definition}\label{def:rectangular_cumulants}
For any monic polynomial $p(x)=x^d+\sum_{i=1}^d{x^{d-i}a_{2i}}$ of degree~$d$ and any $n\in\Z_{\ge 0}$, define its \textbf{$\pmb{(n,d)}$-rectangular (finite free) cumulants} $K^{n,d}_2[p], K^{n,d}_4[p],\dots,K^{n,d}_{2d}[p]$ as the quantities
\begin{equation*}
K^{n,d}_{2\ell}[p] := \ell\cdot[z^{2\ell}]\,\ln\left( 1 + \sum_{i=1}^d{\frac{a_{2i}}{(-d)_i(-d-n)_i}\,z^{2i}} \right),\quad\ell=1,\dots,d.
\end{equation*}
\end{definition}

Alternatively, we can first uniquely define the sequence $\big( K^{n,d}_{2\ell}[p] \big)_{\ell\ge 1}$ by the equality
\begin{equation}\label{eq:alternative}
\exp\left( \sum_{\ell=1}^\infty{ \frac{K^{n,d}_{2\ell}[p]}{\ell} z^{2\ell} } \right)
= 1 + \sum_{i=1}^d{\frac{a_{2i}}{(-d)_i(-d-n)_i}\,z^{2i}},
\end{equation}
and then declare the first $d$ terms of the sequence to be the $(n,d)$-rectangular cumulants of $p(x)$.

Precise formulas for the $(n,d)$-rectangular cumulants $K^{n,d}_{2\ell}[p]$ in terms of the coefficients $a_{2i}$ of $p(x)$, and viceversa, are available: simply replace $t\mapsto -d,\,u\mapsto -d-n$ in Lem.~\ref{lem:a_k}.
As a special case of Lem.~\ref{thm:basic_convolution}, we deduce that the $(n,d)$-rectangular cumulants linearize $\boxplus^n_d$.

\begin{theorem}\label{thm:linearity}
If $p(x), r(x)$ are monic polynomials of degree~$d$ and $n\in\Z_{\ge 0}$, then
\begin{equation*}
K^{n,d}_{2\ell}\big[ p\boxplus_d^n r \big] = K^{n,d}_{2\ell}[p] + K^{n,d}_{2\ell}[r],\quad\text{for all }\ell=1,\dots,d.
\end{equation*}
\end{theorem}

\begin{remark}
In the special case $n=0$, the $(n\!=\!0,\,d)$-rectangular convolution $\boxplus^{n=0}_d$ is called the \emph{asymmetric additive convolution} in~\cite{MSS_2022}.
In this case, our $(n\!=\!0,\,d)$-rectangular cumulants from Def.~\ref{def:rectangular_cumulants} reduce to the \emph{asymmetric cumulants} considered in~\cite{CS_2024}.
\end{remark}

\section{Applications to the asymptotic theory of nonnegative real-rooted polynomials}\label{sec:applications}

\subsection{Limits of (n,d)-rectangular cumulants}

The theorem and corollary in this subsection can be deduced from results on the rectangular finite $R$-transform, obtained in~\cite{G_2024}, but the proofs are different. We recall that for a polynomial $p(x)$ with nonnegative real roots, its symmetric empirical root distribution $\widetilde\mu[p]$ is defined in Eqn.~\eqref{eq:symmetric_distribution}.

\begin{theorem}\label{thm:limit_moments_cumulants}
Let $q\in [1,\infty)$, let $\{p_d(x)\}_{d\ge 1}$ be such that each $p_d(x)$ is a monic polynomial of degree~$d$ with nonnegative real roots, and let $\mu$ be a symmetric probability measure with finite moments of all orders.
Then $\widetilde\mu[p_d]\to\mu$, as $d\to\infty$, in the sense of moments, if and only if
\begin{equation}\label{eq:cumulant_limits}
\lim_{d\to\infty,\, 1+\frac{n}{d}\to q}{ (-d)^{2\ell-1}K^{n,d}_{2\ell}[p_d] } = q^{-\ell}\kappa_{2\ell}^q[\mu],
\quad\textrm{for all $\ell\in\Z_{\ge 1}$}.
\end{equation}
\end{theorem}
\begin{proof}
\textbf{``If direction'': Limits of cumulants implies limits of moments.}
Let $f_d(x)=x^d+\sum_{j=1}^d{a_{2j}^{(d)}x^{d-j}}$ be an arbitrary monic polynomial of degree~$d$ with only nonnegative real roots.
If the even moments of $\widetilde\mu[f_d]$ are denoted $\widetilde{m}_2[f_d],\widetilde{m}_4[f_d],\cdots$, then we claim that
\begin{equation}\label{eq:moments_coeffs_formula}
\exp\left( -d\sum_{k=1}^\infty{ \frac{\widetilde{m}_{2k}[f_d]}{k} z^{2k} } \right) = 1 + \sum_{j=1}^d{ a_{2j}^{(d)}z^{2j} }.
\end{equation}
In fact, if we denote the roots of $f_d(x)$ by $(\alpha_1^{(d)})^2,\dots,(\alpha_d^{(d)})^2$, then $f_d(x)=x^d+\sum_{j=1}^d{a_{2j}^{(d)}x^{d-j}}=\prod_{k=1}^d{\big( x-(\alpha_k^{(d)})^2 \big)}$.
By multiplying both sides by $x^{-d}$ and doing the change of variables $x=z^{-2}$, we get
\begin{align*}
1+\sum_{j=1}^d{a_{2j}^{(d)}z^{2j}} &= \prod_{i=1}^d{\Big( 1 - (\alpha_i^{(d)})^2 z^2 \Big)} = \prod_{i=1}^d{\big( 1 - \alpha_i^{(d)}z \big)\big( 1 + \alpha_i^{(d)}z \big)}\\
&= \exp\left\{ \sum_{i=1}^d\left( \ln\big( 1 - \alpha_i^{(d)}z \big) + \ln\big( 1 + \alpha_i^{(d)}z \big) \right) \right\}\\
&= \exp\left\{ -\sum_{i=1}^d \sum_{\ell=1}^\infty \frac{\Big(\big(\alpha_i^{(d)}\big)^\ell + \big(-\alpha_i^{(d)}\big)^\ell\Big)}{\ell}z^\ell \right\}\\
&= \exp\left\{ -d\sum_{k=1}^\infty \, \frac{\sum_{i=1}^d\frac{1}{2d}\Big(\big(\alpha_i^{(d)}\big)^{2k} + \big(-\alpha_i^{(d)}\big)^{2k}\Big)}{k}z^{2k} \right\},
\end{align*}
which immediately leads to the equality in Eqn.~\eqref{eq:moments_coeffs_formula}; the last equality in the chain of equations above follows because $\big(\alpha_i^{(d)}\big)^\ell + \big(-\alpha_i^{(d)}\big)^\ell=0$, if $\ell$ is odd, so the sum over $\ell\in\Z_{\ge 1}$ can be restricted to the positive even integers by the change of variables $\ell=2k$.

Next, for all $d\in\Z_{\ge 1}$, $n\in\Z_{\ge 0}$, define the sequence $\big(\widetilde{\kappa}_{2\ell}^{n,d}[f_d]\big)_{\ell\ge 1}$ by the equality of formal power series
\begin{equation}\label{eq:modified_cumulants}
\exp\left( \sum_{\ell=1}^\infty{ \frac{\widetilde{\kappa}_{2\ell}^{n,d}[f_d]}{(-d)^{2\ell-1}q^{\ell}} \frac{z^{2\ell}}{\ell} } \right) = 1 + \sum_{j=1}^d{ \frac{a_{2j}^{(d)}}{(-d)_j(-d-n)_j}z^{2j} }.
\end{equation}
Note that equations~\eqref{eq:moments_coeffs_formula}--\eqref{eq:modified_cumulants} are precisely a specialization of \eqref{eqn:a_k}--\eqref{eqn:a_m} where, in particular, $t=-d$ and $u=-d-n$.
Then, by Thm.~\ref{thm:m_k_3}, we have the following formula, for all $k\in\Z_{\ge 1}$:
\begin{multline}\label{eqn:m_k_tilde}
\widetilde{m}_{2k}[f_d] = \sum_{P\in\Luk^\odd(2k)} \prod_{s\ge 0}{(-d-n+s)^{\#\textrm{down steps of $P$ from height }(2s+1)}} \\
\cdot\prod_{s\ge 1}{(-d+s)^{\#\textrm{down steps of $P$ from height }(2s)}}
\prod_{s\ge 1}{\left[ \frac{\widetilde{\kappa}_{2s}^{n,d}[f_d]}{(-d)^{2s-1}q^s} \right]^{\#\textrm{up steps }(1,2s-1)\textrm{ of $P$}}}.
\end{multline}
To simplify~\eqref{eqn:m_k_tilde}, we need some identities on the statistics of odd \L{}ukasiewicz paths. Firstly,
\begin{equation}\label{eq:displacement_zero}
\sum_{s\ge 1}{ (2s-1) \cdot\#\textrm{up steps }(1,2s-1)} = \#\textrm{down steps}
\end{equation}
follows from the fact that any \L{}ukasiewicz path begins and ends at the same height $y=0$ (so the distance traveled going up must equal the distance traveled going down) and odd \L{}ukasiewicz paths only have up steps of the form $(1,2s-1)$, for some $s\in\Z_{\ge 1}$.

Secondly, for any odd \L{}ukasiewicz path of length $2k$, we have
\begin{align}
\sum_{s\ge 1}{ s\cdot\#\textrm{up steps }(1,2s-1)} &= \frac{1}{2}\sum_{s\ge 1}{ (2s)\cdot\#\textrm{up steps }(1,2s-1)}\nonumber\\
&= \frac{1}{2} \left\{ \sum_{s\ge 1}{ (2s-1)\cdot\#\textrm{up steps }(1,2s-1)} + \sum_{s\ge 1}{ \#\textrm{up steps }(1,2s-1)} \right\}\nonumber\\
&= \frac{1}{2}\left\{ \#\textrm{down steps} + \#\textrm{up steps} \right\}\nonumber\\
&= \frac{1}{2}\cdot (2k) = k.\label{eq:length_k}
\end{align}
where the third equality follows from~\eqref{eq:displacement_zero} and the last one follows from the fact that odd \L{}ukasiewicz paths do not have horizontal steps, so the total number of up and down steps of an odd \L{}ukasiewicz path equals its length.

From~\eqref{eq:displacement_zero} and~\eqref{eq:length_k}, we can rewrite equation~\eqref{eqn:m_k_tilde} as
\begin{multline}\label{eqn:m_k_tilde_2}
\widetilde{m}_{2k}[f_d] = q^{-k}\cdot\sum_{P\in\Luk^\odd(2k)} \prod_{s=0}^{k-1}{ \left( 1+\frac{n}{d}-\frac{s}{d} \right)^{\#\textrm{down steps of $P$ from height }(2s+1)} }\\
\cdot\prod_{s=1}^{k-1}{ \left( 1-\frac{s}{d} \right)^{\#\textrm{down steps of $P$ from height }(2s)} }
\prod_{s=1}^{k}{ \big( \widetilde{\kappa}^{n,d}_{2s}[f_d] \big)^{\#\textrm{up steps }(1,2s-1)\textrm{ of $P$}} },
\end{multline}
where we also changed infinite products to finite ones (e.g.~$\prod_{s\ge 0}$ to $\prod_{s=0}^{k-1}$) by using the fact that any $P\in\Luk^\odd(2k)$ has no vertices above height $(2k-1)$ or up steps $(1,2s-1)$ with $s>k$.

Next, for each $k\in\Z_{\ge 1}$, consider the $k$-variate polynomial
\begin{multline}\label{eqn:explicit_pol_1}
P_{2k}^{\,\widetilde{\kappa}\mapsto\widetilde{m}}(x_2,x_4,\dots,x_{2k} \bar d,n,q) :=
q^{-k}\cdot\sum_{P\in\Luk^\odd(2k)} \prod_{s=0}^{k-1}{ \left( 1+\frac{n}{d}-\frac{s}{d} \right)^{\#\textrm{down steps of $P$ from height}\,2s+1} }\\
\cdot\prod_{s=1}^{k-1}{ \left( 1-\frac{s}{d} \right)^{\#\textrm{down steps of $P$ from height}\,2s} }
\prod_{s=1}^{k}{ x_{2s}^{\#\textrm{up steps }(1,2s-1)\textrm{ of $P$}} },
\end{multline}
so that
\begin{equation}\label{eq:pol_1}
\widetilde{m}_{2k}[f_d] = P_{2k}^{\,\widetilde{\kappa}\mapsto\widetilde{m}}
\Big( \widetilde{\kappa}^{n,d}_2[f_d],\dots,\widetilde{\kappa}^{n,d}_{2k}[f_d] \,\Big|\, d,n,q \Big),
\quad\textrm{for all $k\in\Z_{\ge 1}$}.
\end{equation}
Similarly, consider the multivariate polynomials
\begin{equation*}
Q_{2k}^{\,\kappa\mapsto m}(x_2,x_4,\dots,x_{2k}\bar q) = \sum_{\pi\in\NC^\even(2k)}{ q^{-\even(\pi)} \prod_{B\in\pi}{ x_{|B|} } },
\end{equation*}
for all $k\in\Z_{\ge 1}$.
This definition is such that, if we denote the even moments of a symmetric probability measure $\rho$ on $\R$ by $m_{2k}[\rho]$ and its $q$-rectangular free cumulants by $\kappa_{2\ell}^q[\rho]$, then
\begin{equation}\label{eq:pol_3}
m_{2k}[\rho] = Q_{2k}^{\,\kappa\mapsto m}\Big( \kappa_2^q[\rho],\dots,\kappa_{2k}^q[\rho] \,\Big|\, q \Big),
\quad\textrm{for all $k\in\Z_{\ge 1}$}.
\end{equation}

From the definition~\eqref{eqn:explicit_pol_1}, we obtain the following coefficient-wise limit of polynomials
\begin{multline}\label{eqn:m_k_tilde_3}
\lim_{d\to\infty,\, 1+\frac{n}{d}\to q}{ P_{2k}^{\,\widetilde{\kappa}\mapsto\widetilde{m}}(x_2,\dots,x_{2k} \bar d,n,q) } \\
= q^{-k}\cdot\sum_{P\in\Luk^\odd(2k)} { q^{\sum_{s=0}^{k-1}\#\textrm{down steps of $P$ from height }(2s+1)} }
\prod_{s=1}^k{x_{2s}^{\ \#\textrm{up steps }(1,2s-1)\textrm{ of $P$}}}.
\end{multline}
In order to simplify the RHS of~\eqref{eqn:m_k_tilde_3}, we resort to the bijection between $\Luk^\odd(2k)$ and $\NC^\even(2k)$ described in Cor.~\ref{cor:bijection_2}.
If $P\in\Luk^\odd(2k)$ corresponds to $\pi\in\NC^\even(2k)$ under this bijection, then
\begin{align}
\sum_{s=0}^{k-1}\#\,\textrm{down steps of $P$} &\textrm{ from height }(2s+1) = \#\,\textrm{down steps of $P$ from some odd height }\nonumber\\
&= \#\,\textrm{even numbers in $[2k]$ not being the first in their blocks in $\pi$}\nonumber\\
&= k - \#\,\textrm{even numbers in $[2k]$ being the first in their blocks in $\pi$}\nonumber\\
&= k - \even(\pi),\label{eq:bijection_consequence}
\end{align}
where the second equality follows from Cor.~\ref{cor:bijection_2} and the last one is the definition of $\even(\pi)$.
Plugging~\eqref{eq:bijection_consequence} back into~\eqref{eqn:m_k_tilde_3}, we see that the RHS turns precisely into the expression for $Q_{2k}^{\,\kappa\mapsto m}(x_2,\dots,x_{2k}\bar q)$.
Hence, we obtain the coefficient-wise limits of polynomials
\begin{equation}\label{eq:limit_pol_1}
\lim_{d\to\infty,\, 1+\frac{n}{d}\to q}{ P_{2k}^{\,\widetilde{\kappa}\mapsto\widetilde{m}}(x_2,\dots,x_{2k}\bar d,n,q) }
= Q_{2k}^{\,\kappa\mapsto m}(x_2,\dots,x_{2k}\bar q),\qquad\textrm{for all $k\in\Z_{\ge 1}$}.
\end{equation}

We can now conclude the argument.
By comparing~\eqref{eq:modified_cumulants} with Def.~\ref{def:rectangular_cumulants}, we have
\begin{equation}\label{eq:widekappa_K}
\widetilde{\kappa}^{n,d}_{2\ell}[p_d] = (-d)^{2\ell-1}q^\ell\cdot K^{n,d}_{2\ell}[p_d],\quad\textrm{for all }\ell=1,2,\dots,d.
\end{equation}
Assume that the limits~\eqref{eq:cumulant_limits} hold; by~\eqref{eq:widekappa_K}, this can be stated as
\begin{equation}\label{eq:limit_cumulants}
\lim_{d\to\infty,\,1+\frac{n}{d}\to q}{ \widetilde{\kappa}^{n,d}_{2\ell}[p_d] } = \kappa_{2\ell}^q[\mu],\quad\textrm{for all }\ell\in\Z_{\ge 1}.
\end{equation}
Hence, by Eqn.~\eqref{eq:pol_1} for the polynomials $p_d$, Eqn.~\eqref{eq:pol_3} for $\mu$, and the coefficient-wise limit of polynomials~\eqref{eq:limit_pol_1}, we obtain the limits of even moments
\begin{equation}\label{eq:limit_moments_0}
\lim_{d\to\infty}{ \widetilde{m}_{2k}[p_d] } = m_{2k}[\mu],\quad\textrm{for all }k\in\Z_{\ge 1}.
\end{equation}
This suffices to prove the desired limit $\widetilde\mu[p_d]\to\mu$, in the sense of moments, since all these measures are symmetric and have vanishing odd moments.

\medskip
\noindent\textbf{``Only if direction'': Limits of moments implies limits of cumulants.}
For this part, assuming the limits~\eqref{eq:limit_moments_0}, we want to conclude~\eqref{eq:limit_cumulants}.
First, note that each $P_{2k}^{\,\widetilde{\kappa}\mapsto\widetilde{m}}(x_2,\dots,x_{2k}\bar d,n,q)$ is of degree $2k$, if we set $\deg x_{2s}:=2s$, for $s=1,\dots,k$; indeed, this follows from~\eqref{eq:length_k}.
This means that $P_{2k}^{\,\widetilde{\kappa}\mapsto\widetilde{m}}(x_2,\dots,x_{2k}\bar d,n,q)$ is a linear combination of terms $x_\pi=\prod_{B\in\pi}{x_{|B|}}$, as $\pi$ ranges over $\P^\even(2k)$.
Moreover, the system of equations~\eqref{eq:pol_1} (for $f_d\mapsto p_d$) can be enlarged to a system where each $\widetilde{m}_{\sigma}[p_d] := \prod_{B\in\sigma}{\widetilde{m}_{|B|}[p_d]}$, for $\sigma\in\P^\even(2k)$, is expressed as a linear combination of terms $\widetilde{\kappa}^{n,d}_\pi[p_d] := \prod_{B\in\pi}{\widetilde{\kappa}^{n,d}_{|B|}[p_d]}$, for $\pi\in\P^\even(2k)$.
In matrix terms, this can be written as
\begin{equation}\label{eq:matrix_1}
\vec{\mathbf{m}}[p_d] = \mathbf{P}(d,n,q)\cdot\vec{\pmb{\kappa}}^{\,n,d}[p_d],
\end{equation}
where $\vec{\mathbf{m}}[p_d]$ and $\vec{\pmb{\kappa}}^{\,n,d}[p_d]$ are column vectors of size $|\P^\even(2k)|$, with entries $\kappa^{n,d}_\sigma[p_d]$ and $m_\pi[p_d]$, respectively.
The matrix $\mathbf{P}(d,n,q)$ is of size $|\P^\even(2k)|\times|\P^\even(2k)|$.
Due to the fact that each polynomial $P_{2k}^{\,\widetilde{\kappa}\mapsto\widetilde{m}}(x_2,\dots,x_{2k} \bar d,n,q)$ has the form
\begin{equation}\label{eq:triangular_form}
\begin{aligned}
P_{2k}^{\,\widetilde{\kappa}\mapsto\widetilde{m}}(x_2,\dots,x_{2k} \bar d,n,q)
=&\ q^{-k}\cdot\prod_{s=0}^{k-1}\left(1 + \frac{n}{d} - \frac{s}{d}\right)
\cdot\prod_{s=1}^{k-1}\left(1 - \frac{s}{d}\right)\cdot x_{2k}\\
&+ \textrm{some polynomial in the variables }x_2,\dots,x_{2k-2},
\end{aligned}
\end{equation}
we deduce that $\mathbf{P}(d,n,q)$ is lower triangular (if we order rows and columns with the reverse lexicographic order) and it has nonzero diagonal entries.
The equation~\eqref{eq:matrix_1} can thus be inverted:
\begin{equation*}
\vec{\pmb{\kappa}}^{\,n,d}[p_d] = \mathbf{P}(d,n,q)^{-1}\cdot\vec{\mathbf{m}}[p_d].
\end{equation*}
By looking at the last entry (corresponding to the set partition with only one block of size $2k$), we see $\widetilde{\kappa}^{n,d}_{2k}[f_d]$ on the left and a linear combination of $\widetilde{m}_\pi[f_d]$, $\pi\in\P^\even(2k)$, on the right.
In other words, the last row of $\mathbf{P}(d,n,q)^{-1}$ furnishes a polynomial $P_{2k}^{\,\widetilde{m}\mapsto\widetilde{\kappa}}(x_2,x_4,\dots,x_{2\ell}\bar d,n,q)$ of degree $2k$ such that
\begin{equation}\label{eq:pol_2}
\widetilde{\kappa}^{n,d}_{2k}[p_d] = P_{2k}^{\,\widetilde{m}\mapsto\widetilde{\kappa}}
\Big( \widetilde{m}_2[p_d],\dots,\widetilde{m}_{2k}[p_d] \,\Big|\, d,n,q \Big).
\end{equation}
As $k$ was arbitrary, these polynomials exist for all $k\in\Z_{\ge 1}$.

Similarly, from the definition, note that each $Q_{2k}^{\,\kappa\mapsto m}(x_2,\dots,x_{2k}\bar q)$ is of degree $2k$, if we set $\deg x_{2s}=2s$.
Then, as before, we obtain a matrix equality
\begin{equation}\label{eq:system_2}
\vec{\mathbf{m}}[\mu] = \mathbf{Q}(q)\cdot\vec{\pmb{\kappa}}^{\,q}[\mu]
\end{equation}
where $\vec{\mathbf{m}}[\mu]$, $\vec{\pmb{\kappa}}^{\,q}[\mu]$ are column vectors of size $|\P^\even(2k)|$ and entries $m_\pi[p_d]$, $\kappa^{n,d}_\sigma[p_d]$, while $\mathbf{Q}(q)$ is a matrix of size $|\P^\even(2k)|\times|\P^\even(2k)|$.
Since each $Q_{2k}^{\,\kappa\mapsto m}$ is of the form
\begin{equation}\label{eq:triangular_form_2}
Q_{2k}^{\,\kappa\mapsto m}(x_2,\dots,x_{2k} \bar d,n,q) = x_{2k} + \textrm{some polynomial in the variables }x_2,\dots,x_{2k-2},
\end{equation}
then $\mathbf{Q}(q)$ is an lower uni-triangular matrix, in particular, it is invertible and the system~\eqref{eq:system_2} can be inverted:
\begin{equation*}
\vec{\pmb{\kappa}}^{\,q}[\mu] = \mathbf{Q}(q)^{-1}\cdot\vec{\mathbf{m}}[\mu]
\end{equation*}
Looking at the last entry of this equality, we see $\kappa^q[\mu]$ on the left and a linear combination of terms $m_\pi[\mu]=\prod_{B\in\pi}{m_{|B|}[\mu]}$, $\pi\in\P^\even(2k)$, on the right.
As a result, there exists a polynomial $Q_{2k}^{\,m\mapsto\kappa}(x_2,\dots,x_{2k}\bar q)$ of degree $2k$ such that
\begin{equation}\label{eq:pol_4}
\kappa_{2k}^q[\mu] = Q_{2k}^{\,m\mapsto\kappa}\Big( m_2[\mu],\dots,m_{2k}[\mu] \,\Big|\, q \Big),
\end{equation}
and since $k$ was arbitrary, this is true for all $k\in\Z_{\ge 1}$.

Next, by~\eqref{eq:limit_pol_1} and the construction of the matrices $\mathbf{P}(d,n,q)$, $\mathbf{Q}(q)$, we have $\lim{\mathbf{P}(d,n,q)} = \mathbf{Q}(q)$ entry-wise, where the limit is taken as $d\to\infty$ and $1+\frac{n}{d}\to q$.
Observe that the coefficient of $x_{2k}$ in~\eqref{eq:triangular_form} converges to $1$ (in the same limit regime), so each diagonal entry of the upper triangular matrix $\mathbf{P}(d,n,q)$ converges to $1$, implying $\lim{\det\mathbf{P}(d,n,q)} = 1$.
By Cramer's rule, it follows that $\lim{\mathbf{P}(d,n,q)^{-1}} = \mathbf{Q}(q)^{-1}$ entry-wise.
In particular, by looking at the last row of this matrix equality, we deduce the coefficient-wise limits of polynomials:
\begin{equation}\label{eq:limit_pol_2}
\lim_{d\to\infty,\, 1+\frac{n}{d}\to q}{ P_{2k}^{\,\widetilde{m}\mapsto\widetilde{\kappa}}(x_2,\dots,x_{2k} \bar d,n,q) }
= Q_{2k}^{\,m\mapsto\kappa}(x_2,\dots,x_{2k} \bar q),\quad\textrm{for all $k\in\Z_{\ge 1}$}.
\end{equation}
Finally, \eqref{eq:pol_2}, \eqref{eq:pol_4} and \eqref{eq:limit_pol_2} readily show the desired fact that the limits of moments~\eqref{eq:limit_moments_0} imply the limits of cumulants~\eqref{eq:limit_cumulants}.
The proof is finished.
\end{proof}

\begin{corollary}\label{cor:d_infinity}
Let $q\in[1,\infty)$, let $\mu,\nu$ be compactly supported symmetric probability measures on $\R$, and let $\{p_d(x)\}_{d\ge 1}$, $\{r_d(x)\}_{d\ge 1}$ be such that $p_d(x)$, $r_d(x)$ are monic polynomials of degree~$d$ with all their roots being real and nonnegative, for all $d\in\Z_{\ge 1}$.
Further, assume that we have the limits $\widetilde\mu[p_d]\to\mu$, $\widetilde\mu[r_d]\to\nu$, as $d\to\infty$, in the sense of moments.
Then $\widetilde\mu\big[ p_d\boxplus^n_d r_d \big]\to\mu\boxplus_q\nu$, in the sense of moments, in the regime where $n,d\to\infty$ and $1+\frac{n}{d}\to q$.
\end{corollary}
\begin{proof}
The proof here continues the proof of Thm.~\ref{thm:limit_moments_cumulants}.
The assumption that $\widetilde\mu[p_d]\to\mu$, $\widetilde\mu[r_d]\to\nu$, in the sense of moments, implies:
\begin{equation}\label{eq:limit_moments_2}
\lim_{d\to\infty}{ \widetilde{m}_{2k}[p_d] } = m_{2k}[\mu],\qquad
\lim_{d\to\infty}{ \widetilde{m}_{2k}[r_d] } = m_{2k}[\nu],\qquad
\textrm{for all }k\in\Z_{\ge 1}.
\end{equation}
As a result, by~\eqref{eq:pol_2}, \eqref{eq:pol_4} and the coefficient-wise limit of polynomials~\eqref{eq:limit_pol_2}, the limits~\eqref{eq:limit_moments_2} imply
\begin{equation}\label{eq:limit_cumulants_2}
\lim_{d\to\infty,\, 1+\frac{n}{d}\to q}{ \widetilde{\kappa}^{n,d}_{2\ell}[p_d] } = \kappa^q_{2\ell}[\mu],\qquad
\lim_{d\to\infty,\, 1+\frac{n}{d}\to q}{ \widetilde{\kappa}^{n,d}_{2\ell}[r_d] } = \kappa^q_{2\ell}[\nu],\qquad
\textrm{for all }\ell\in\Z_{\ge 1}.
\end{equation}
Next, for all $\ell\in\Z_{\ge 1}$, we have
\begin{multline}\label{eq:limits_cumulants}
\lim_{d\to\infty,\, 1+\frac{n}{d}\to q}{ \widetilde{\kappa}^{n,d}_{2\ell}\big[p_d \boxplus^n_d r_d\big] }
= \lim_{d\to\infty,\, 1+\frac{n}{d}\to q}{ \Big( \widetilde{\kappa}^{n,d}_{2\ell}[p_d]  + \widetilde{\kappa}^{n,d}_{2\ell}[r_d] \Big) }\\
= \kappa_{2\ell}^q[\mu] + \kappa_{2\ell}^q[\nu] = \kappa_{2\ell}^q[\mu\boxplus_q\nu],
\end{multline}
where the second equality is a consequence of~\eqref{eq:limit_cumulants_2}, the third one is the result of our ad-hoc Theorem-Definition~\ref{thm_def}, and the first one (recalling from~\eqref{eq:widekappa_K} that $\widetilde{\kappa}^{n,d}_{2\ell}[p_d] = (-d)^{2\ell-1}q^\ell\cdot K^{n,d}_{2\ell}[p_d]$, for all $\ell=1,\dots,d$) follows from Thm.~\ref{thm:linearity}.
Finally, by~\eqref{eq:pol_1}, \eqref{eq:pol_3} and the coefficient-wise limit of polynomials~\eqref{eq:limit_pol_1}, the limits~\eqref{eq:limits_cumulants} imply
\begin{equation*}
\lim_{d\to\infty,\,1+\frac{n}{d}\to q}{ \widetilde{m}_{2k}\big[p_d \boxplus^n_d r_d\big] } = m_{2k}[\mu\boxplus_q\nu],\quad
\textrm{for all }k\in\Z_{\ge 1},
\end{equation*}
which proves the desired convergence $\widetilde\mu\big[p_d \boxplus^n_d r_d\big]\to\mu\boxplus_q\nu$, in the sense of moments.
\end{proof}

\subsection{Proof of the Main Result~III}\label{sec:asymptotic_real_roots}

We need some preparations for the proof.
Let $\R[x]$ be the space of real polynomials on $x$ and denote the derivative by $D=\frac{d}{dx}$.
If $n\in\Z_{\ge 0}$, note that the operator $x^{-n}Dx^{n+1}D$ acts on $\R[x]$ and decreases the degree by $1$.
In fact, $(x^{-n}Dx^{n+1}D)x^i = i(i+n)x^{i-1}$, and more generally,
\begin{equation}\label{eq:rectangular_operator}
(x^{-n}Dx^{n+1}D)^k x^i = (-i)_k (-i-n)_k\cdot x^{i-k},\quad\textrm{for all }k=0,1,\dots,i.
\end{equation}
For any polynomial $p(x)=\sum_{k=0}^d{a_{2k}x^{d-k}}$ of degree~$d$, there exists another polynomial $P(x)$ of the same degree such that $p(x)=P(x^{-n}Dx^{n+1}D)x^d$. In fact, by~\eqref{eq:rectangular_operator}, it follows that
\begin{equation}\label{eq:pol_P}
P(x) = \sum_{k=0}^d{ \frac{a_{2k}x^k}{(-d)_k(-d-n)_k} }\ 
\Longrightarrow\ p(x)=P(x^{-n}Dx^{n+1}D)x^d.
\end{equation}
If $r(x)$ is another polynomial of degree~$d$ and $R(x)$ is the polynomial of the same degree such that $r(x)=R(x^{-n}Dx^{n+1}D)x^d$, then by~\eqref{eq:rectangular_convolution}, \eqref{eq:rectangular_operator} and \eqref{eq:pol_P}, it is easy to verify that
\begin{equation}\label{eq:observation}
\begin{aligned}
(p\boxplus^n_dr)(x) &= P(x^{-n}Dx^{n+1}D)R(x^{-n}Dx^{n+1}D)x^d\\
&= P(x^{-n}Dx^{n+1}D)r(x) = R(x^{-n}Dx^{n+1}D)p(x).
\end{aligned}
\end{equation}

\begin{proof}[\textbf{Proof of the Main Result~III}]
Since $x^{-n}Dx^{n+1}D$ decreases the degree of any polynomial by $1$, then $\exp\Big(-\frac{s^2}{n} x^{-n}Dx^{n+1}D\Big)$ is a valid operator on $\R[x]$ and acts like
\begin{equation}\label{eq:application_vertex}
\exp\left(-\frac{s^2}{n} x^{-n}Dx^{n+1}D\right)p_d(x)
= \sum_{k=0}^d{ \frac{(-s^2)^k}{n^k k!} \big( x^{-n}Dx^{n+1}D \big)^k p_d(x) },
\end{equation}
on the degree~$d$ polynomial $p_d(x)$. The right hand side of~\eqref{eq:application_vertex} can be rewritten in terms of the polynomial
\begin{equation*}
R_{n,d}(x) := \sum_{k=0}^d{ \frac{ (-s^2)^k }{n^k k!}\,x^k }
\end{equation*}
as $R_{n,d}(x^{-n}Dx^{n+1}D)p_d(x)$. Hence, by~\eqref{eq:observation}, if we define $r_{n,d}(x) := R_{n,d}(x^{-n}Dx^{n+1}D)x^d$, then
\begin{equation*}
\exp\left(-\frac{s^2}{n} x^{-n}Dx^{n+1}D\right)p_d(x) = \big( p_d \boxplus^n_d r_{n,d} \big)(x).
\end{equation*}
By virtue of Thm.~\ref{thm:linearity}, we then have
\begin{equation}\label{eq:additive_cumulants}
K^{n,d}_{2\ell}\left[ \exp\left(-\frac{s^2}{n} x^{-n}Dx^{n+1}D\right)p_d \right] = K^{n,d}_{2\ell}[p_d] + K^{n,d}_{2\ell}[r_{n,d}],
\quad\textrm{for all }\ell=1,\dots,d.
\end{equation}
By using~\eqref{eq:rectangular_operator}, the polynomial $r_{n,d}(x)$ can be explicitly computed as:
\begin{equation}\label{r_polynomial}
r_{n,d}(x) = R_{n,d}(x^{-n}Dx^{n+1}D)x^d = \sum_{k=0}^d{ \frac{(-s^2)^k (-d)_k(-d-n)_k}{n^k k!}\, x^{d-k} }.
\end{equation}
Thus, by Def.~\ref{def:rectangular_cumulants}, the $(n,d)$-rectangular cumulants $K^{n,d}_{2\ell}[r_{n,d}]$ are
\begin{equation*}
K^{n,d}_{2\ell}[r_{n,d}] = \ell\cdot[z^{2\ell}]\ln\left( 1 + \sum_{k=1}^d{ \frac{(-s^2)^k}{n^k k!}z^{2k} } \right),
\quad\textrm{for }\ell=1,\dots,d.
\end{equation*}
For any $\ell\le d$, note that
\begin{multline*}
[z^{2\ell}]\ln\left( 1 + \sum_{k=1}^d{ \frac{(-s^2)^k}{n^k k!} z^{2k} } \right)
= [z^{2\ell}]\ln\left( 1 + \sum_{k=1}^\infty{ \frac{(-s^2)^k}{n^k k!} z^{2k} } \right)\\
= [z^{2\ell}]\ln\left(\exp(-s^2z^2/n)\right) = [z^{2\ell}](-s^2z^2/n) = -\delta_{\ell, 1}\cdot\frac{s^2}{n},
\end{multline*}
therefore $\displaystyle K^{n,d}_{2\ell}[r_{n,d}] = -\delta_{\ell, 1}\cdot\frac{s^2}{n}$, for all $\ell=1,\dots,d$.
Plugging this back into~\eqref{eq:additive_cumulants} gives
\begin{equation}\label{eq:additive_cumulants_2}
K^{n,d}_{2\ell}\left[ \exp\left(-\frac{s^2}{n} x^{-n}Dx^{n+1}D\right)p_d \right]
= K^{n,d}_{2\ell}[p_d] - \delta_{\ell, 1}\cdot\frac{s^2}{n},\quad\textrm{for all }\ell=1,\dots,d.
\end{equation}
By the assumption of the theorem and the ``only if direction'' of Thm.~\ref{thm:limit_moments_cumulants}, we have
\begin{equation*}
\lim_{d\to\infty,\, 1+\frac{n}{d}\to q}{ (-d)^{2\ell-1}K^{n,d}_{2\ell}[p_d] } = q^{-\ell}\kappa_{2\ell}^q[\mu],
\end{equation*}
for all $\ell\in\Z_{\ge 1}$. Then by~\eqref{eq:additive_cumulants_2}, we have the limits
\begin{multline}\label{eq:final_1}
\lim_{d\to\infty,\, 1+\frac{n}{d}\to q}{ (-d)^{2\ell-1} K^{n,d}_{2\ell}\left[ \exp\left(-\frac{s^2}{n} x^{-n}Dx^{n+1}D\right)p_d \right] }\\
= \lim_{d\to\infty,\, 1+\frac{n}{d}\to q}{ (-d)^{2\ell-1} \left( K^{n,d}_{2\ell}[p_d] - \delta_{\ell, 1}\cdot\frac{s^2}{n} \right) }
= q^{-\ell}\kappa_{2\ell}^q[\mu] + \delta_{\ell, 1}\cdot\frac{s^2}{q-1},
\end{multline}
for all $\ell\in\Z_{\ge 1}$.
As the symmetric probability measure $\lambda_{qs^2/(q-1)}^{(q)}$ from Lem.~\ref{lem:helpful_measure} has $q$-rectangular free cumulants $\kappa^q_{2\ell}\Big[ \lambda_{qs^2/(q-1)}^{(q)} \Big] = \delta_{\ell, 1}\cdot\frac{qs^2}{q-1}$, we can write
\begin{equation}\label{eq:final_2}
\delta_{\ell, 1}\cdot\frac{s^2}{q-1} = q^{-\ell}\cdot\delta_{\ell, 1}\cdot\frac{qs^2}{q-1} = q^{-\ell}\cdot\kappa^q_{2\ell}\Big[ \lambda_{qs^2/(q-1)}^{(q)} \Big].
\end{equation}
Consequently,
\begin{align*}
\lim_{d\to\infty,\, 1+\frac{n}{d}\to q}{ (-d)^{2\ell-1} K^{n,d}_{2\ell}\left[ \exp\left(-\frac{s^2}{n} x^{-n}Dx^{n+1}D\right)p_d \right] }
&= q^{-\ell} \Big( \kappa_{2\ell}^q[\mu] + \kappa^q_{2\ell}\Big[ \lambda_{qs^2/(q-1)}^{(q)} \Big] \Big)\\
&= q^{-\ell}\kappa_{2\ell}^q \Big[ \mu \boxplus_q \lambda_{qs^2/(q-1)}^{(q)} \Big],
\end{align*}
where the first equality follows from~\eqref{eq:final_1}--\eqref{eq:final_2} and the second from Theorem-Definition~\ref{thm_def}.
Finally, the ``if direction'' of Thm.~\ref{thm:limit_moments_cumulants} yields the desired limit~\eqref{eq:conclusion_3} in the sense of moments.
\end{proof}

\begin{remark}
The theory of the rectangular finite free convolution is related to the classical Laguerre orthogonal polynomials.
Indeed, recall that the Laguerre polynomial $L_d^{(\alpha)}(x)$ with parameter $\alpha>-1$ and degree $d$ is
\[
L_d^{(\alpha)}(x) = \sum_{k=0}^d{ {d+\alpha \choose d-k} \frac{(-x)^k}{k!}}.
\]
Then the polynomial $r_{n,d}(x)$ in Eqn.~\eqref{r_polynomial} inside the proof above can be expressed in terms of the Laguerre polynomial $L_d^{(n)}(x)$ with parameter $n\in\Z_{\ge 0}$ by the equality
\[
r_{n,d}(x) = \left(-\frac{s^2}{n}\right)^d\cdot d!\cdot L_d^{(n)}\left(\frac{nx}{s^2}\right).
\]
Moreover, the unique monic polynomial of degree $d$ with $(n,d)$-rectangular cumulants being $K_{2\ell}^{n,d}=\delta_{\ell,1}$ is exactly $d!\cdot L_d^{(n)}(-x)$, as can be verified from~\eqref{eq:alternative}.
This is then the limiting object for the following rectangular finite free version of the Central Limit Theorem:

\emph{If $p(x)$ is any monic polynomial of degree $d$ with second cumulant $K_2^{n,d}$ equal to $1$, then the polynomials
\[
N^{-d}\cdot \big(p\boxplus_d^n\cdots\boxplus_d^n p\big)(Nx) \quad\text{(N-fold convolution)}
\]
converge coefficient-wise, as $N\to\infty$, to the polynomial $d!\cdot L_d^{(n)}(-x)$.}

The symmetric finite free version of this result was stated and proved in \cite[Thm.~6.7]{M_2021}, and then revisited in \cite[Example~6.1]{AP_2018}, with the use of finite free cumulants.
The rectangular version stated above can be proved by exactly the same argument as in the second reference, so we omit the details.
\end{remark}

\section{q-generalizations}\label{sec:q_analogues}

In this section,\footnote{Since $q$ has already appeared as the parameter for the theory of rectangular free probability, we use the bold~$\q$ in this section to denote the ``quantum'' deformation parameter.} we consider a $\q$-analogue of the Main Result~I from the introduction; in order to state it, we need some classical notations from $\q$-analysis, e.g. from~\cite{GR_2011}.
The $\q$-numbers are $[n]_\q=(1-\q^n)/(1-\q)$, $n\in\Z_{\ge 1}$, while $[0]_\q=1$; the $\q$-factorial is $n!_\q=[1]_\q[2]_\q\cdots[n]_\q$, for $n\in\Z_{\ge 1}$, while $0!_\q=1$; the $q$-exponential function is $\exp_\q(x) = \sum_{n=0}^\infty{x^n/n!_\q}$.

The $\q$-derivative of a function $f(x)$ is $D_\q f(x)=\frac{f(x)-f(\q x)}{(1-\q)x}$; for example, $D_\q\exp_\q(x)=\exp_\q(x)$.
The $\q$-Leibniz rule is $D_\q\big(f(x)g(x)\big) = f(x)\cdot D_\q g(x) + D_\q f(x)\cdot g(\q x)$. On the other hand, any naive chain rule for $D_\q$ fails.
However, Gessel~\cite{G_1982} showed that there exists a chain rule under a proper $\q$-analogue of function composition; we next follow the presentation from~\cite{J_1996}.
Given a function $f(x)$ with $f(0)=0$, we can write $f(x)=\sum_{n=1}^\infty{ f_n\frac{x^n}{n!_\q} }$, for some values $f_1, f_2,\cdots$. The \emph{$\q$-symbolic powers} of $f(x)$ are defined inductively by $f^{\q;[0]}(x):=1$ and
\begin{equation*}
D_\q f^{\q;[k]}(x):=[k]_\q\cdot f^{\q;[k-1]}(x)D_\q f(x),\quad f^{\q;[k]}(0):=0,\quad\textrm{for all }k\in\Z_{\ge 1}.
\end{equation*}
For example, $x^{\q;[k]}=x^k$ and more generally
\begin{equation}\label{eqn:symbolic_powers}
f^{\q;[k]}(x) = x(1-\q^k)\sum_{n=0}^\infty{\q^n\cdot f^{\q;[k-1]}(x\q^n)\cdot (D_\q f)(x\q^n)},\quad\text{for all }k\in\Z_{\ge 1},
\end{equation}
by induction on $k$.
Then, if $f(x)=\sum_{n=1}^\infty{ f_n\frac{x^n}{n!_\q} }$, $g(x)=\sum_{n=1}^\infty{ g_n\frac{x^n}{n!_\q} }$, the \emph{$\q$-composition} of $f$ and $g$ is defined as
\begin{equation}\label{eqn:composition}
g[f]:=\sum_{n=1}^\infty{g_n\frac{f^{\q;[n]}}{n!_\q}}.
\end{equation}
Note that the $\q$-symbolic powers and $\q$-composition depend on $\q$ (though the $\q$-dependence of the $\q$-composition is not shown in its notation).
With these definitions, the chain rule is: $D_\q\big(g[f]\big) = (D_\q g)[f]\cdot D_\q f$.
In particular, for $g=\exp_\q$, we have
\begin{equation}\label{eqn:q_chain}
D_\q\big(\exp_\q[f]\big) = \big(\exp_\q[f]\big)\cdot D_\q f.
\end{equation}
To state our theorem, recall the $\q$-Pochhammer symbol $(t;\q)_n:=\prod_{i=0}^{n-1}(1-t\q^i)$, for $n\in\Z_{\ge 1}$.

\begin{theorem}[$\q$-analogue of Main Result~I]\label{thm:q_analogue}
Let $\q,t,u$ be formal parameters and let $(a_2,a_4,\dots)$, $(\kappa_2,\kappa_4,\dots)$, $(m_2,m_4,\dots)$ be sequences that are related to one another by:
\begin{align}
\exp_{\q^{-\frac{1}{2}}}\!\left[ \sum_{\ell=1}^\infty{ \frac{\kappa_{2\ell}}{[\ell]_\q} \q^{\ell-\frac{1}{2}} z^{2\ell} } \right]
&= 1 + \sum_{n=1}^\infty{ a_{2n}\frac{(1-\q)^{2n}}{(t;\q)_n(u;\q)_n} \q^{n^2-\frac{3n}{2}} z^{2n} },\label{q_equations_1}\\
\exp_{\q^{-\frac{1}{2}}}\!\left[ \frac{1-t}{1-\q}\sum_{k=1}^\infty{ \frac{m_{2k}}{[k]_\q} z^{2k} } \right] &= 1 + \sum_{n=1}^\infty{ a_{2n}z^{2n} }.\label{q_equations_2}
\end{align}
In the relations above, we used the $\q$-numbers $[\ell]_\q=\frac{1-\q^\ell}{1-\q}$, $[k]_\q=\frac{1-\q^k}{1-\q}$, but $\q^{-\frac{1}{2}}$-composition (namely, definitions \eqref{eqn:symbolic_powers}--\eqref{eqn:composition} should be employed with the parameter $\q^{-\frac{1}{2}}$ instead of $\q$) for the expressions of the form $\exp_{\q^{-\frac{1}{2}}}[f(z)]$ on the left hand sides of \eqref{q_equations_1}--\eqref{q_equations_2}.

\begin{enumerate}[label=(\alph*)]
	\item For all $k\in\Z_{\ge 1}$, we have
\begin{equation}\label{eq:q_ops}
m_{2k} = \q^{\frac{1}{2}}\cdot [z^0]\Big( \q\,\Delta_{\q;t,u} T_{\q^{-\frac{1}{2}}} + \q^{\frac{1}{2}}*_gT_{\q^{-\frac{1}{2}}} \Big)^{2k-1}\!\left(g(z)\right),
\end{equation}
where $g(z):=\sum_{\ell=1}^\infty{\kappa_{2\ell} z^{2\ell-1}}$; $T_{\q^{-\frac{1}{2}}}$ is the shift operator $T_{\q^{-\frac{1}{2}}}f(z):=f\big(\q^{-\frac{1}{2}}z\big)$; $[z^0]f(z)$ is the constant term of the power series $f(z)\in\R[[z]]$; and finally, $\Delta_{\q;t,u}$ is the linear operator on $\R[[z]]$ defined uniquely by its action on monomials:
\begin{equation}\label{eq:Delta}
\begin{gathered}
\Delta_{\q;t,u}z^{2n+1} := \frac{1 - u\q^n}{1-\q}z^{2n},\quad n\in\Z_{\ge 0},\\
\Delta_{\q;t,u}z^{2n} := \frac{1 - t\q^n}{1-\q}z^{2n-1},\quad n\in\Z_{\ge 1},\qquad \Delta_{\q;t,u}(1) := 0.
\end{gathered}
\end{equation}

	\item For all $k\in\Z_{\ge 1}$, we have
\begin{multline}\label{eq:q_moment}
m_{2k} = \sum_{P\in\Luk^\odd(2k)} 
\q^{\#\textrm{down steps of $P$} + \frac{1}{2}(\#\textrm{up steps of $P$} - \textrm{sum of heights of vertices of $P$})}\\
\cdot\prod_{s\ge 0}{\Big((1-u\q^s)/(1-\q)\Big)^{\#\textrm{down steps of $P$ from height }(2s+1)}}\\
\cdot\prod_{s\ge 1}{\Big((1-t\q^s)/(1-\q)\Big)^{\#\textrm{down steps of $P$ from height }(2s)}}
\prod_{s\ge 1}{\kappa_{2s}^{\#\textrm{up steps }(1,2s-1)\textrm{ of $P$}}}.
\end{multline}
\end{enumerate}
For formula~\eqref{eq:q_moment}, we used the terminology that the vertex $(x,y)$ of $P\in\Luk^\odd(2k)$ has \textbf{height $y$} and the \textbf{sum of heights of vertices of $P$} is the sum of all $y$-coordinates of vertices of $P$.
\end{theorem}
\begin{proof}[Sketch of proof]
For part (a), repeat the same proof as for Thm.~\ref{thm:m_k_2}, but with the operators $d_{t,u}$, $\partial_{t,u}$ on $\R[[z]]$ used in that proof replaced by $d_{\q;t,u}$, $\partial_{\q;t,u}$, defined by
\begin{gather*}
d_{\q;t,u}z^{2n+1} := \frac{1-u\q^{-\frac{1}{2}}}{\q^{-1}-1}z^{2n},\quad n\in\Z_{\ge 0},\\
d_{\q;t,u}z^{2n} := \frac{1-t}{\q^{-1}-1}z^{2n-1},\quad n\in\Z_{\ge 1},\qquad d_{\q;t,u}(1) := 0,
\end{gather*}
and $\displaystyle\partial_{\q;t,u} := \frac{1}{\q^{-\frac{1}{2}}+1}D_{\q^{-\frac{1}{2}}} + d_{\q;t,u}$, or explicitly,
\begin{equation}\label{eq:q_partial}
\begin{gathered}
\partial_{\q;t,u} z^{2n+1} := \frac{\q^{-n-\frac{1}{2}}-u\q^{-\frac{1}{2}}}{\q^{-1}-1}z^{2n},\quad n\in\Z_{\ge 0},\\
\partial_{\q;t,u} z^{2n} := \frac{\q^{-n}-t}{\q^{-1}-1}z^{2n-1},\quad n\in\Z_{\ge 1},\qquad \partial_{\q;t,u}(1) := 0.
\end{gathered}
\end{equation}
We only comment that the $\q$-version of the chain rule~\eqref{eqn:q_chain} is essential in the argument and that the necessary analogue of Claim~\ref{claim:step_3} is the relation
\begin{multline*}
\partial_{\q;t,u}\Big( h(z)\exp_{\q^{-\frac{1}{2}}}[G(z)] \Big) \\
= \Big( \partial_{\q;t,u} + \q^{\frac{1}{2}}*_gT_{\q^{-\frac{1}{2}}} \Big)(h(z))\cdot\exp_{\q^{-\frac{1}{2}}}[G(z)] + [z^0]h(z)\cdot d_{\q;t,u}\Big(\exp_{\q^{-\frac{1}{2}}}[G(z)]\Big),
\end{multline*}
where $G(z)=\sum_{\ell=1}^\infty{ \frac{\kappa_{2\ell}}{[\ell]_\q} \q^{\ell-\frac{1}{2}} z^{2\ell} }$, which leads to
\begin{equation*}
m_{2k} = \q^{\frac{1}{2}}\cdot[z^0] \Big( \partial_{\q;t,u} + \q^{\frac{1}{2}}*_g T_{\q^{-\frac{1}{2}}} \Big)^{2k-1}\big(g(z)\big),
\end{equation*}
for all $k\in\Z_{\ge 1}$.
Finally, the desired equation~\eqref{eq:q_ops} follows because $\partial_{\q;t,u}$ and $\q\,\Delta_{\q;t,u}T_{\q^{-\frac{1}{2}}}$ are equivalent operators, as seen from formulas \eqref{eq:Delta} and \eqref{eq:q_partial}.

Part (b) follows from part (a) by a similar reasoning as in the proof of Thm.~\ref{thm:m_k_3}.
One new argument needed is that the factor $\q$ in the operator $\q\,\Delta_{\q;t,u}T_{\q^{-\frac{1}{2}}}$ in the RHS of \eqref{eq:q_ops} accounts for the factor $\q^{\#\textrm{down steps}}$ in equation~\eqref{eq:q_moment}, whereas the factor $\q^{\frac{1}{2}}$ in $\q^{\frac{1}{2}}*_gT_{\q^{-\frac{1}{2}}}$ accounts for the factor $\q^{\frac{1}{2}\#\textrm{up steps}}$ in~\eqref{eq:q_moment}, and finally, the $\q^{-\frac{1}{2}}$-shift $T_{\q^{-\frac{1}{2}}}$ in both $\q\,\Delta_{\q;t,u}T_{\q^{-\frac{1}{2}}}$ and $\q^{\frac{1}{2}}*_gT_{\q^{-\frac{1}{2}}}$ leads to the factor $\q^{-\frac{1}{2}\textrm{sum of heights of vertices}}$, after the proper combinatorial interpretation.
\end{proof}

\begin{example}
If $\kappa_n = \delta_{n,2}\cdot\q^{-\frac{1}{2}}c$, then
\begin{equation*}
G(z) = \sum_{\ell=1}^\infty{ \frac{\kappa_{2\ell}}{[\ell]_\q} \q^{\ell-\frac{1}{2}} z^{2\ell} } = cz^2,\qquad g(z) = \sum_{\ell=1}^\infty{\kappa_{2\ell}z^{2\ell-1}} = \q^{-\frac{1}{2}}cz.
\end{equation*}
A simple induction shows
\begin{equation*}
G^{\,\q^{-\frac{1}{2}}; [n]}(z) = \frac{\q^{\frac{1}{2}n(n-1)}\cdot n!_{\q^{-\frac{1}{2}}}}{n!_\q}\,c^nz^{2n},\quad \text{for all }n\in\Z_{\ge 1},
\end{equation*}
and therefore
\begin{equation*}
\exp_{\q^{-\frac{1}{2}}}[G(z)] = \sum_{n=0}^\infty{ \frac{G^{\,\q^{-\frac{1}{2}}; [n]}}{n!_{\q^{-\frac{1}{2}}}} }
= 1 + \sum_{n=1}^\infty{ \frac{\q^{\frac{1}{2}n(n-1)}}{n!_\q}\,c^nz^{2n} }.
\end{equation*}
The left hand side above is the same as the left hand side in Eqn.~\eqref{q_equations_1}, therefore the right hand sides must also agree; as a result, by comparing the coefficients of $z^{2n}$, we find that
\begin{equation}\label{eq:a_2n}
a_{2n} = \frac{\q^{n-\frac{n^2}{2}}(t;\q)_n(u;\q)_n}{(1-\q)^{2n}\,n!_\q}\,c^n,\quad\textrm{for all } n\in\Z_{\ge 1}.
\end{equation}
By taking the $\q^{-\frac{1}{2}}$-derivative $D_{\q^{-\frac{1}{2}}}$ to both sides of~\eqref{q_equations_2} and applying the chain rule~\eqref{eqn:q_chain}, we obtain
\[
\left(1 + \sum_{i=1}^\infty{ a_{2i}z^{2i} }\right)\cdot\left( \frac{1-t}{1-\q}\sum_{k=1}^\infty{ \frac{m_{2k}}{[k]_\q} D_{\q^{-\frac{1}{2}}} z^{2k} } \right) = \sum_{n=1}^\infty{ a_{2n}\, D_{\q^{-\frac{1}{2}}}z^{2n} }
\]
By using $D_{\q^{-\frac{1}{2}}} z^{2k} = \q^{-k+\frac{1}{2}}[k]_{\q}\cdot (1+\q^{\frac{1}{2}})$ and $D_{\q^{-\frac{1}{2}}} z^{2n} = \q^{-n+\frac{1}{2}}[n]_{\q}\cdot (1+\q^{\frac{1}{2}})$, it follows that
\[
\frac{1-t}{1-\q} \left(1 + \sum_{i=1}^\infty{ a_{2i}z^{2i} }\right)\left(\sum_{k=1}^\infty{ m_{2k}\,\q^{-k+\frac{1}{2}} z^{2k-1} } \right) = \sum_{n=1}^\infty{ a_{2n}\,\q^{-n+\frac{1}{2}}[n]_{\q}\, z^{2n-1} }.
\]
By comparing the coefficients of $z^{2n-1}$ on both sides, we deduce
\begin{equation}\label{eq:a_2n_2}
\frac{1-t}{1-\q}\left\{ m_{2n} + \sum_{\ell=1}^{n-1}{ \q^{n-\ell}m_{2\ell}\,a_{2n-2\ell} } \right\} = [n]_\q\cdot a_{2n},
\quad\textrm{for all }n\in\Z_{\ge 1}.
\end{equation}
With \eqref{eq:a_2n} and \eqref{eq:a_2n_2}, we can then inductively find $m_{2n}$; for example,
\begin{equation}\label{first_ms}
m_2 = \q^{\frac{1}{2}}\frac{1-u}{1-\q}c,\qquad
m_4 = \frac{(1-u)(1+\q-\q t-\q u)}{(1-\q)^2} c^2.
\end{equation}
We shall verify that $m_2,m_4$ match the formulas obtained by their combinatorial interpretation in~\eqref{eq:q_moment}.
First, observe that there exists only one \L{}ukasiewicz path in $\Luk^\odd(2)$, namely the one in Fig.~\ref{fig:luk1}, and it has weight
\begin{equation*}
\q^{1+\frac{1}{2}(1-1)}\frac{1-u}{1-\q}\kappa_2 = \q^{\frac{1}{2}}\frac{1-u}{1-\q}c,
\end{equation*}
which matches exactly the value of $m_2$ in~\eqref{first_ms}.
\begin{figure}[h]
\begin{center}
\includegraphics[width=0.19\textwidth]{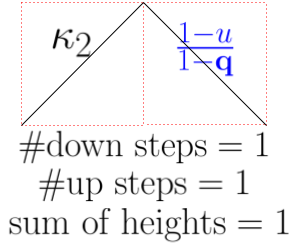}
\caption{The only \L{}ukasiewicz path in $P\in\Luk^\odd(2)$.}\label{fig:luk1}
\end{center}
\end{figure}

Next, since $\kappa_n=0$, whenever $n\ne 2$, the only \L{}ukasiewicz paths that contribute to the formula~\eqref{eq:q_moment} will be the ones that only have up steps of size $1$ and there are exactly two of them in $\Luk^\odd(4)$ with that property, namely the ones depicted in Fig.~\ref{fig:luk2}.
They have weights
\begin{figure}[h]
\begin{center}
\includegraphics[width=0.52\textwidth]{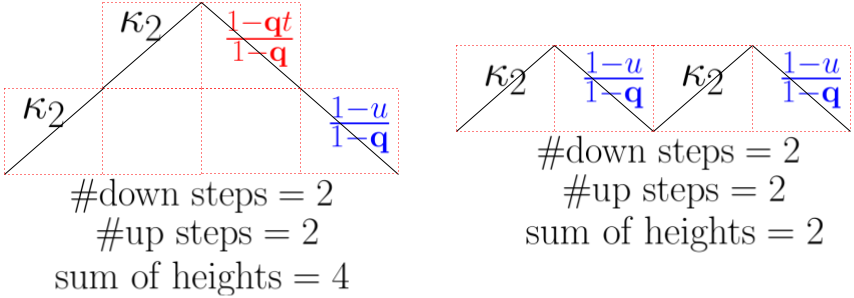}
\caption{The \L{}ukasiewicz paths in $\Luk^\odd(4)$ that only have up steps of size $1$.}\label{fig:luk2}
\end{center}
\end{figure}
\begin{equation*}
\q^{2+\frac{1}{2}(2-4)}\frac{(1-\q t)(1-u)}{(1-\q)^2}\kappa_2^2 = \frac{(1-\q t)(1-u)}{(1-\q)^2}c^2,\qquad
\q^{2+\frac{1}{2}(2-2)}\frac{(1-u)^2}{(1-\q)^2}\kappa_2^2 = \q\frac{(1-u)^2}{(1-\q)^2}c^2,
\end{equation*}
which indeed add up to the value of $m_4$ in~\eqref{first_ms}.
\end{example}

\section{Appendix: beta-deformed singular values and eigenvalues}\label{appendix}

We briefly explain the occurrence of Eqns.~\eqref{eqn:a_k}--\eqref{eqn:a_m} in the topic of the Appendix's title.

\subsection{High temperature beta-singular values}\label{subsec:appendix_1}

The \emph{$q$-rectangular free convolution} was defined by Benaych-Georges~\cite{BG_2009} in the following setting.
For a matrix $A\in\C^{M\times N}$, $M\le N$, with singular values $a_1\ge\dots\ge a_M\ge 0$, its symmetric empirical measure is defined as
\begin{equation*}
\widetilde\mu_A := \frac{1}{2M}\sum_{i=1}^M{(\delta_{a_i} + \delta_{-a_i})}.
\end{equation*}
Let $\{A_M\}_{M\ge 1}$, $\{B_M\}_{M\ge 1}$ be sequences of independent $M\times N(M)$ random matrices, with invariant laws under the natural $U(M)\times U(N(M))$-actions, with deterministic singular values, and such that $N(M)\to\infty$, $N(M)/M\to q\in[1,\infty)$, as $M\to\infty$. Assume that $\widetilde\mu_{A_M}\to\mu_A$, $\widetilde\mu_{B_M}\to\mu_B$, weakly as $M\to\infty$, for two compactly supported symmetric probability measures $\mu_A, \mu_B$ on $\R$.
Then~\cite{BG_2009} shows that the symmetric empirical measures $\widetilde\mu_{C_M}$ of the sums $C_M=A_M+B_M$ converge weakly, in probability, to another probability measure, uniquely determined by $\mu_A$ and $\mu_B$, and he called it the $q$-rectangular free convolution of $\mu_A$ and $\mu_B$:
\begin{equation}\label{eqn:rectangular_limit}
\widetilde\mu_{C_M}\to\mu_A\boxplus_q\mu_B,\quad\textrm{as }M\to\infty.
\end{equation}

In~\cite{X_2023}, the author found a one-parameter deformation to the operation of $q$-rectangular free convolution (so there are two deformation parameters $q,\gamma$ in total) by considering a novel map that takes as input two $M$-tuples $\a=(a_1\ge\dots\ge a_M\ge 0)$, $\b=(b_1\ge\dots\ge b_M\ge 0)$, as well as parameters $\beta>0$ (the ``inverse temperature''), $N\in\Z$, $N>M$, and outputs a generalized function on $M$-tuples $\cc=(c_1\ge\dots\ge c_M\ge 0)$.
The map was regarded by the author as a ``$\beta$-generalization of addition of rectangular matrices'', because when $\beta=2$, it can be interpreted as what happens when considering singular values of sums of independent $M\times N$ random matrices.
To be more precise, the map mentioned originates by regarding the function $\C^{N\times M}\ni X\mapsto\E\big[ \exp(\Re\,Tr(AX)) \big]$ as a Fourier-transform of $A$, since for the sum $C=A+B$ of independent matrices, one has
\begin{equation*}
\E\big[ \exp(\Re\,\textrm{Tr}(CX)) \big] = \E\big[ \exp(\Re\,\textrm{Tr}(AX)) \big]\cdot\E\big[ \exp(\Re\,\textrm{Tr}(BX)) \big].
\end{equation*}
Note that $\E\big[ \exp(\Re\,\textrm{Tr}(AX)) \big]=:\phi_{M,N}(\a,\x)$ depends only on the singular values $\a=(a_1,\dots,a_M)$, $\x=(x_1,\dots,x_M)$ of $A,X$, and the dimension $N$.
The \emph{multivariate BC-type Bessel function} $\phi^{(\beta)}_{M,N}(\a,\x)$ furnishes a $\beta$-generalization, with a matrix interpretation only for $\beta=1,2,4$, that serves to define, for any $\beta>0$, the \emph{$\beta$-sum of $M\times N$ random matrices} $\a\boxplus^{\beta}_{M,N}\!\b$, as the unique generalized function $\E$ on $\R^M$ such that
\begin{equation*}
\E\Big[ \phi^{(\beta)}_{M,N}\Big( \a\boxplus^\beta_{M,N}\!\b,\x \Big) \Big] 
= \phi^{(\beta)}_{M,N}(\a,\x)\cdot\phi^{(\beta)}_{M,N}(\b,\x),\ \text{ for all }\x\in\R^M.
\end{equation*}
In the high temperature limit regime
\begin{equation*}
\beta\to 0,\qquad
\frac{\beta M}{2}\to\gamma\in (0,\infty),\qquad
\frac{N}{M}\to q\in[1,\infty),
\end{equation*}
\cite{X_2023} proved a LLN that generalizes the previous limit~\eqref{eqn:rectangular_limit}, where the new limiting object depends on a two-parameter binary operation that the author called the \emph{$(q,\gamma)$-convolution} of probability measures, and denoted it by $\mu_\a\boxplus_{q,\gamma}\mu_\b$.

If we consider the moment sequences $m^\a=(m_1^\a,m_2^\a,\dots)$, $m^\b=(m_1^\b,m_2^\b,\dots)$, $m^\cc=(m_1^\cc, m_2^\cc,\dots)$ of $\mu_\a$, $\mu_\b$ and $\mu_\cc:=\mu_\a\boxplus_{q,\gamma}\mu_\b$, respectively, then the $(q,\gamma)$-convolution map $(m^\a,m^\b)\mapsto m^\cc$ turns out to be equivalent, after the moment-to-coefficient transformation~\eqref{eqn:a_m}, to the convolution of sequences~\eqref{eqn:convolution_intro}, upon the specialization of parameters:
\begin{equation}\label{eq:tu_1}
t=\gamma,\qquad u=q\gamma.
\end{equation}
Moreover, the coefficient-to-cumulant transformation~\eqref{eqn:a_k} is what led~\cite{X_2023} to define the $(q,\gamma)$-cumulants, which linearize the $(q,\gamma)$-convolution of measures.

\subsection{High temperature beta-eigenvalues}\label{subsec:appendix_2}

In~\cite{BGCG_2022}, the authors studied a $\beta$-dependent map that takes as input $\a=(a_1,\dots,a_N)$, $\b=(b_1,\dots,b_N)$ and outputs a generalized function on $\R^N$, denoted $\a\boxplus^\beta_N\b$.
When $\beta=1,2,4$, this coincides with the sum of independent orbital-distributed self-adjoint $N\times N$ random matrices with fixed eigenvalues, while for general $\beta>0$, it is defined by means of Fourier-type transforms and multivariate Bessel functions of type~A.
The LLN proved in~\cite{BGCG_2022} shows that in the regime $\beta\to 0$, $N\to\infty$, $\frac{\beta N}{2}\to\gamma\in (0,\infty)$, the operation $\boxplus^\beta_N$ tends to the \emph{$\gamma$-convolution} $\boxplus_\gamma$ of probability measures. It turns out that if $\mu_\cc=\mu_\a\boxplus_\gamma\mu_\b$, and we let $\m_\a:=(m^\a_1,m^\a_2,\dots)$, $\m_\b=(m_1^\b,m^\b_2,\dots)$, $\m_\cc=(m^\cc_1,m^\cc_2,\dots)$ be their moment sequences, then the map $(\m^\a,\m^\b)\mapsto\m_\cc$ coincides with the convolution of sequences~\eqref{eqn:convolution_intro} upon the transformation given by~\eqref{def_cm} and the specialization of parameters:
\begin{equation}\label{eq:tu_2}
t=\gamma,\qquad u=0.
\end{equation}
The cumulants given the formulas~\eqref{eqn:a_k}--\eqref{eqn:a_m} are exactly the \emph{$\gamma$-cumulants} introduced in~\cite{BGCG_2022}.
The setting in this section is similar to that of Sec.~\ref{subsec:appendix_1}, but for square Hermitian (not rectangular) matrices.
In terms of formulas, note that by comparing \eqref{eq:tu_1}--\eqref{eq:tu_2}, the setting here is a strict specialization of the previous section (when $q=0$), yet even in this case the moment-cumulant formulas found in~\cite{BGCG_2022} could only be proved by means of Dunkl theory.

One last worthwhile comment is that the generating function identities \eqref{eqn:a_k}--\eqref{eqn:a_m} (for $u=0$) were first derived in~\cite{BGCG_2022} by taking limits in the high temperature regime of the integral representations of specialized multivariate Bessel functions, see~\cite[Thms.~5.1--5.2]{C_2021}.
It is plausible that the analogous integral representations for Jack polynomials~\cite[Thms.~2.2--2.3]{C_2018b} and Macdonald polynomials~\cite[Thm.~3.2--3.3]{C_2018a} are related to our $\q$-analogue of Sec.~\ref{sec:q_analogues}.


\begin{thebibliography}{9}

\bibitem[AP-18]{AP_2018}
Octavio Arizmendi and Daniel Perales. Cumulants for finite free convolution. Journal of Combinatorial Theory, Series~A~155 (2018), pp.~244--266.

\bibitem[AGVP-23]{AGVP_2023}
Octavio Arizmendi, Jorge Garza-Vargas and Daniel Perales. Finite free cumulants: Multiplicative convolutions, genus expansion and infinitesimal distributions. Transactions of the American Mathematical Society~376, no.~06 (2023), pp.~4383--4420.

\bibitem[BG-07]{BG_2007}
Florent Benaych-Georges. Infinitely divisible distributions for rectangular free convolution: classification and matricial interpretation. Probability Theory and Related Fields~139 (2007), pp.~143--189.

\bibitem[BG-09]{BG_2009}
Florent Benaych-Georges. Rectangular random matrices, related convolution. Probability Theory and Related Fields~144 (2009), pp.~471--515.

\bibitem[BGCG-22]{BGCG_2022}
Florent Benaych-Georges, Cesar Cuenca and Vadim Gorin. Matrix addition and the Dunkl transform at high temperature. Communications of Mathematical Physics~394 (2022), pp.~735--795.

\bibitem[CS-24]{CS_2024}
Volodymyr Chub and Stanislav Surmylo. Cumulants for asymmetric additive convolution. Final report, Yulia's Dream Program, MIT PRIMES, Cambridge MA (2024). Available at: \url{https://drive.google.com/file/d/1qxRdEwgZKf9rjL-d-vNu9d9MjbWy3y8P/view}.

\bibitem[Cue-18a]{C_2018a}
Cesar Cuenca. Asymptotic Formulas for Macdonald Polynomials and the Boundary of the $(q,t)$-Gelfand-Tsetlin graph. SIGMA~14 (2018), 001.

\bibitem[Cue-18b]{C_2018b}
Cesar Cuenca. Pieri Integral Formula and Asymptotics of Jack Unitary Characters. Selecta Mathematica, New Series, vol.~24, no.~3 (2018), pp. 2737--2789.

\bibitem[Cue-21]{C_2021}
Cesar Cuenca. Universal Behavior of the Corners of Orbital Beta Processes. International Mathematics Research Notices~2021, no.~19, pp.~14761--14813.

\bibitem[CDM-23]{CDM_2023}
Cesar Cuenca, Maciej Do\l{}\k{e}ga and Alexander Moll. Universality of global asymptotics of Jack-deformed random Young diagrams at varying temperatures. Annals of Probability, Vol.~54, No.~1 (2026), pp.~421--488.

\bibitem[CD-25a]{CD_2025a}
Cesar Cuenca and Maciej Do\l{}\k{e}ga. Discrete $N$-particle systems at high temperature through Jack generating functions Preprint; arXiv:2502.13098 (2025).

\bibitem[CD-25b]{CD_2025b}
Cesar Cuenca and Maciej Do\l{}\k{e}ga. Crystallization of discrete $N$-particle systems at high temperature. Preprint; arXiv:2510.23496 (2025).

\bibitem[GR-11]{GR_2011}
George Gasper and Mizan Rahman. Basic hypergeometric series. Vol.~96. Cambridge University Press, 2011.

\bibitem[Ges-82]{G_1982}
Ira M.~Gessel. A $q$-analog of the exponential formula. Discrete Mathematics~40, no.~1 (1982), pp.~69--80.

\bibitem[Gri-22]{G_2022}
Aurelien X.~Gribinski. Rectangular finite free probability theory. PhD thesis, Princeton University, 2022.

\bibitem[Gri-24]{G_2024}
Aurelien Gribinski. A Theory of Singular Values for Finite Free Probability. Journal of Theoretical Probability~37, no.~2 (2024), pp.~1257--1298.

\bibitem[GM-22]{GM_2022}
Aurelien Gribinski and Adam W.~Marcus. A rectangular additive convolution for polynomials. Combinatorial Theory~2, no.~1 (2022).

\bibitem[Joh-96]{J_1996}
Warren P.~Johnson. Some applications of the $q$-exponential formula. Discrete Mathematics~157, no.~1 (1996), pp.~207--225.

\bibitem[Kab-22]{K_2022}
Zakhar Kabluchko. Lee-Yang zeroes of the Curie-Weiss ferromagnet, unitary Hermite polynomials, and the backward heat flow. Annales Henri Lebesgue~8 (2025), pp.~1--34.

\bibitem[KS-09]{KS_2009}
Claus K\"{o}stler and Roland Speicher. A noncommutative de Finetti theorem: Invariance under quantum permutations is equivalent to freeness with amalgamation. Communications in Mathematical Physics~291, no.~2 (2009), pp.~473--490.

\bibitem[LT-16]{LT_2016}
Francois Lemeux and Pierre Tarrago. Free wreath product quantum groups: the monoidal category, approximation properties and free probability. Journal of Functional Analysis~270, no.~10 (2016), pp.~3828--3883.

\bibitem[Mar-21]{M_2021}
Adam W.~Marcus. Polynomial convolutions and (finite) free probability. Preprint; arXiv:2108.07054 (2021).

\bibitem[MSS-22]{MSS_2022}
Adam W.~Marcus, Daniel A.~Spielman and Nikhil Srivastava. Finite free convolutions of polynomials. Probability Theory and Related Fields~182, no.~3-4 (2022), pp.~807--848.

\bibitem[MS-17]{MS_2017}
James A.~Mingo and Roland Speicher. Free probability and random matrices. Vol.~35. New York: Springer, 2017.

\bibitem[Mir-21]{Mir_2021}
Benjamin B.~P.~Mirabelli. Hermitian, non-hermitian and multivariate finite free probability. PhD thesis, Princeton University, 2021.

\bibitem[NS-16]{NS_2006}
Alexandru Nica and Roland Speicher. London Mathematical Society Lecture Note Series~335. Cambridge University Press, Cambridge (2006).

\bibitem[Rai-05]{R_2005}
Eric M.~Rains. BC$_n$-symmetric polynomials. Transformation Groups~10 (2005), pp.~63--132.

\bibitem[Spe-94]{S_1994}
Roland Speicher. Multiplicative functions on the lattice of non-crossing partitions and free convolution. Mathematische Annalen~298, no.~1 (1994), pp.~611--628.

\bibitem[Sta-12]{S_2012}
Richard P.~Stanley. Enumerative Combinatorics, Vol.~I. Cambridge Studies in Advanced Mathematics~49, Cambridge University Press, Cambridge, second edition (2012).

\bibitem[Voi-91]{V_1991}
Dan Voiculescu. Limit laws for random matrices and free products. Inventiones Mathematicae~104, no.~1 (1991), pp.~201--220.

\bibitem[Xu-23]{X_2023}
Jiaming Xu. Rectangular matrix additions in low and high temperatures. Preprint; arXiv:2303.13812 (2023).

\end{thebibliography}
\end{document}